\documentclass{article}
\usepackage{fullpage}
\usepackage{graphicx}
\usepackage{dirtytalk}
\usepackage{amsmath,amssymb,bm}
\usepackage{xcolor}
\usepackage[colorlinks,linkcolor=blue,citecolor=magenta]{hyperref}
\usepackage{doi}
\usepackage[nameinlink]{cleveref}% After hyperref

\newcommand{\minimize}{\operatornamewithlimits{minimize}}

\newcommand{\subg}{\bm{\mathrm{g}}}  % bold roman g (for vector)
\newcommand{\xb}{\bm{\mathrm{x}}}  % bold roman x (for vector)
\newcommand{\yb}{\bm{\mathrm{y}}}  % bold roman x (for vector)
\newcommand{\eb}{\bm{\mathrm{e}}}  % bold roman x (for vector)
\newcommand{\x}{x}  % plain x (for components)
\newcommand{\nq}{n_q} % number of quadrupoles
\newcommand{\nd}{n_d} % number of dipoles
\newcommand{\nac}{n_{\rm ac}} % number of algebraic constraints
\newcommand{\nsc}{n_{\rm sc}} % number of simulation-based constraints
\newcommand{\bb}{\mathbf{b}} % bold b
\newcommand{\cb}{\mathbf{c}} % bold c
\newcommand{\db}{\mathbf{d}} % bold d
\newcommand{\Ab}{\mathbf{A}} % bold A
\newcommand{\Yb}{\mathbf{Y}} % bold Y
\newcommand{\eps}{\varepsilon} 
\newcommand{\epsb}{\bm{\eps}} % bold eps
\newcommand{\R}{\mathbb{R}} % Reals 
\newcommand{\zerob}{\bm{0}}     % Vector of zeros 
\newcommand{\lambdab}{\bm{\lambda}}     % bold lambda
\newcommand{\gammab}{\bm{\gamma}}     % bold gamma
\newcommand{\betab}{\bm{\beta}}     % bold beta
\newcommand{\nx}{n_{\bf{x}}} % Number of decision variables

\usepackage{authblk}

\begin{document}

\title{Derivative-Free Optimization of a Rapid-Cycling Synchrotron\thanks{This  manuscript  is  based  upon  work  supported  by  the  applied  mathematics and the Scientific Discovery through Advanced Computing (SciDAC) programs of  the  Office  of  Advanced  Scientific  Computing  Research,  Office  of  Science,  U.S.\ Department of Energy, under Contract DE-AC02-06CH11357. This manuscript has been authored by Fermi Research Alliance, LLC under Contract  DE-AC02-07CH11359 with the U.S.\ Department of Energy, Office of Science, Office of High Energy Physics. Synergia development has been supported by the Office of Advanced Scientific Computing Research and Office of High Energy Physics SciDAC program.}
}


\author[1]{Jeffrey S.~Eldred}
\author[2]{Jeffrey Larson}
\author[2,3]{Misha Padidar}
\author[1]{Eric Stern}
\author[2]{Stefan M.~Wild}
\affil[1]{Fermi National Accelerator Laboratory, Batavia, IL \linebreak[4] \texttt{\small jseldred@fnal.gov}; \texttt{\small egstern@fnal.gov}}
\affil[2]{Mathematics \& Computer Science Division, Argonne National Laboratory, Lemont, IL \linebreak[4] \texttt{\small jmlarson@anl.gov}; \texttt{\small wild@anl.gov}}
\affil[3]{Center for Applied Mathematics, Cornell University, Ithaca, NY \linebreak[4] \texttt{\small map454@cornell.edu}}

% \author{Jeffrey S.~Eldred    \and %https://orcid.org/0000-0003-4432-072X
%         Jeffrey Larson        \and % ORCID:0000-0001-9924-2082
%         Misha Padidar      \and
%         %https://orcid.org/0000-0002-0710-4377
%         Eric Stern \and %https://orcid.org/0000-0003-4317-7159
%         Stefan~M.~Wild}
%          %SW ORCID:0000-0002-6099-2772


% \institute{J.S.\ Eldred \and  E. Stern \at
%               %P.O.\ Box 0500, 60510
%               Fermi National Accelerator Laboratory, Batavia, IL, USA \\
%               \email{jseldred@fnal.gov}; \email{egstern@fnal.gov}           %  \\
% %             \emph{Present address:} of F. Author  %  if needed
%            \and
%             J.\ Larson \and M.\ Padidar \and S.M.\ Wild \at
%               Mathematics \& Computer Science Division, Argonne National Laboratory, Lemont, IL, USA \\
%               \email{jmlarson@anl.gov}; \email{wild@anl.gov}           %  \\
%            \and
%             M.\ Padidar \at
%               Center for Applied Mathematics, Cornell University, Ithaca, NY, USA\\
%               \email{map454@cornell.edu}           %  \\             
% }

\date{\today}

\maketitle

\begin{abstract}
We develop and solve a constrained optimization model to identify an integrable optics rapid-cycling synchrotron lattice design that performs well in several capacities. Our model encodes the design criteria into 78 linear and nonlinear constraints, as well as a single nonsmooth objective, where the objective and some constraints are defined from the output of Synergia, an accelerator simulator. We detail the difficulties of the 23-dimensional simulation-constrained decision space and establish that the space is nonempty.  We use a derivative-free manifold sampling algorithm to account for structured nondifferentiability in the objective function. Our numerical results quantify the dependence of solutions on constraint parameters and the effect of the form of objective function. 

%GP  - just curious: wouldn't one expect the  solutions  to depend  on  the  constraint  parameters?
%SW: - I changed illustrate to quantify, but if the constraint is inactive, then small changes in the parameters would not change the solution.
% \keywords{Numerical optimization \and Simulation optimization \and Particle accelerator design \and Rapid cycling synchrotron}
\end{abstract}

\section{Introduction}
\label{sec:intro}

With the advent of the Long-Baseline Neutrino Facility (LBNF) at Fermilab~\cite{DUNE}, there is a strong motivation to follow the PIP-II 1.2~MW upgrade~\cite{PIP2} with a 2.4~MW upgrade of the Fermilab proton accelerator complex. The construction of a new high-intensity rapid-cycling synchrotron (RCS) accelerator would provide a clear path to achieve the 2.4~MW benchmark for the LBNF program and set the stage for the next generation of particle physics experiments~\cite{EldredSyphers,EldredJINST,Nagaitsev}.

Designing such an RCS accelerator can be greatly aided by numerical optimization techniques, which have 
been used to
design and tune particle accelerators. Particle accelerator optics and performance characteristics have been optimized by simplex methods~\cite{Huang2018,Borland2005}; particle swarms~\cite{HuangSafranek2014,Pang2014}; multiobjective genetic algorithms~\cite{Sun2017,Yang2011}; and, more recently, genetic algorithms enhanced by machine learning~\cite{Edelen2020,Li2018}. Particle accelerator operations also deploy online tuning methods, which include model-based methods~\cite{MintyBook} and local extrema-seeking methods~\cite{Neveu2019,Neveu2017,Scheinker2014} and methods employing Gaussian process models~\cite{Duris2020,Roussel2021}. In this paper we develop a nonsmooth, constrained optimization model for particle accelerator performance and present the first nonsmooth optimization of nonlinear integrable accelerator lattice optics.

%Due to the recent success of numerical optimization in designing and tuning particle accelerators \cite{Neveu2019,Neveu2017,Li2018,Edelen2020,Yang2011,Huang2018,Borland2005,Pang2014,Huang2014,Scheinker2014,Duris2020,Roussel2021},
We develop and solve an optimization model that seeks an RCS lattice design that performs well in several capacities. Our model takes the form
\begin{equation}
    \label{eq:genprob}
    \begin{array}{l}
    \displaystyle    \minimize_{\xb\in \Omega} f(\xb)\\
    \Omega = \left\{\xb\in \R^{\nx} : \Ab \xb \leq \bb, \; c_i(\xb) \leq 0, i =1,\ldots, \nsc\right\}.
        \end{array}
\end{equation}
The objective $f:\R^{\nx}\mapsto\R$ and nonlinear constraints $\cb:\R^{\nx}\mapsto\R^{\nsc}$ encode desired lattice properties through the simulation-based evaluation of a lattice defined by the $\nx$ decision variables $\xb$, which represent the lengths, strengths, and positions of the lattice elements: quadrupoles, dipoles, radio-frequency cavity (RF) inserts, and nonlinear (NL) inserts. Describing this problem with an optimization model provides a geometrically meaningful formulation of the decision space and an understandable parametric description for the tradeoffs between desired lattice properties. However, with the benefits of an interpretable problem description come difficulties in finding points that are feasible for the problem. This manuscript seeks to serve as a guide for future RCS lattice designs by describing the construction and solution of this simulation-based optimization model as well as obstacles in solving the optimization problem.

We describe a principled search of the feasible set $\Omega$. 
We are especially careful to distinguish constraint functions that are algebraically available from those that depend on simulation output. 
The use of Synergia, an accelerator modeling framework, is key in this study because it allows for rapid simulation of a lattice design through linearized optics. From this basic framework of linear accelerator optics, particle accelerator performance can be assessed theoretically with additional considerations for adiabatic changes, machine errors, coupling effects, nonlinearities, parametric resonances, particle interactions, feedback systems, and operational requirements. For example, Wei~\cite{Wei2003} gives an overview of design considerations for intense hadron synchrotron accelerators.

As a case study, we develop and solve an optimization model for an RCS compatible with nonlinear integrable optics. Nonlinear integrable optics is a recent particle accelerator technology that enables strong nonlinear focusing without generating new parametric resonances~\cite{Danilov}. A promising application of integrable optics is in high-intensity rings, where nonlinearities are known to suppress the formation of beam halos~\cite{EldredJINST,EldredIPAC17,Webb2012} and enhance Landau damping of charge-dominated collective instabilities~\cite{Macridin:2015vua}. Over the next several years, the application of nonlinear integrable optics for intense beams will be studied experimentally at the Fermilab Accelerator Science and Technology Integrable Optics Test Accelerator (IOTA)~\cite{Antipov2017,Valishev2020}, as well as smaller-scale test accelerators at the University of Maryland~\cite{Ruisard2019} and at the Rutherford Appleton Laboratory~\cite{Martin2019}.

This paper is structured as follows. \Cref{sec:rcs} details the RCS lattice design problem and describes the particular synchrotron features desired in the case study considered in this paper. \Cref{sec:synergia} details the use of Synergia~\cite{Synergia,Amundson:2004qd} to simulate these features. \Cref{sec:study} formulates a mathematical model of the decision space $\Omega$. 
Through algebraic reductions, we are able to reduce our optimization model to a minimal set of 55 linear constraints and 23 simulation-based constraints involving 32 decision variables. 
\Cref{sec:feasibility} and \Cref{sec:objective} respectively  address two of the key challenges in solving \cref{eq:genprob}: finding a feasible design in $\Omega$ and tackling nondifferentiability in objective functions of interest.  In both cases we take advantage of known dependence of the constraint and objective functions on various Synergia outputs, which allows us to find and certify locally optimal solutions.  \Cref{sec:results} shows the results of our numerical investigation and explores the dependence of our solutions on problem parameters. \Cref{sec:conc} concludes the paper with a look to future work.

\section{Rapid-cycling synchrotron design case study}
\label{sec:rcs}

For the Fermilab 2.4~MW RCS upgrade, using an RCS lattice design compatible with nonlinear integrable optics has been considered in order to enhance the capabilities of the proton facility, reduce technical risk, and dramatically reduce costs associated with constructing the superconducting linac injector~\cite{Eldred2020,EldredJINST}. It is also important to demonstrate that the additional constraints associated with nonlinear integrable optics do not adversely affect conventional lattice design criteria.

An overview of fundamental particle accelerator dynamics and notation is given in~\cite{SyphersBook} and~\cite{LeeBook}. Particle accelerators maintain the stability of particle beams by providing focusing in all three planes of motion. Longitudinal focusing is provided by RF resonating cavities, which accelerate particles in the beam selectively depending on their arrival time. In the transverse planes, steering is provided by dipole magnets, and linear (de)focusing is provided by quadrupole magnets.

The linear transverse ``optics'' of the particle accelerator is described by Hill's equations. The piecewise explicit time-dependence of Hill's equations represent the passage of the particle beam through the sequence of focusing magnets known as the accelerator lattice.  Hill's equations are solved in their canonical form by Floquet's theorem, which gives rise to optics functions including the amplitude-like beta function and phase advance $\phi$ for each plane of motion. The transverse particle motion is known as a ``betatron oscillation,'' and the number of oscillations in one revolution around a particle accelerator ring is known as the ``betatron tune.'' Particle dynamics in an accelerator can also be represented by a sequence of transfer matrices, which can be thought of as mappings between snapshots of solutions to Hill's equations at specific locations.
%GP  -  Hill's equations  or equation?
%SW: I'm going with Hill's equations based on the web search: "hill's equations" particle accelerator
% SW: I made this consistent throughout
The optics functions can be derived as properties of these matrices. Physical and engineering constraints, such as the length and strength of components, can be related to the linear accelerator optics only through the interaction of many independent accelerator magnets; consequently, nonlinear optimization methods are relied on to explore the decision space.

In application, the transverse optics are defined for a particle at a reference momentum (which can change as the beam accelerates), and additional terminology is needed to describe the relative motion of particles at other momentums. The particle trajectory that requires no restorative focusing force is known as the beam ``orbit,'' and the ``dispersion'' is defined to be the change in beam orbit with respect to beam momentum. The change in the betatron tune with beam momentum is known as ``chromaticity.'' The ``momentum compactor factor'' describes the change in the path length with momentum, aggregated over one revolution around a particle accelerator ring.

In~\cite{Danilov} the ``Danilov--Nagaitsev'' criteria for integrable accelerator design first require an alternating sequence of linear and nonlinear sections. The linear sections, referred to as T inserts, are a sequence of dipole and quadrupole magnets with an overall $\pi$-integer betatron phase advance in the horizontal and vertical plane. In the nonlinear sections, the lattice should be dispersion-free, and the horizontal and vertical beta functions should be equal to each other throughout. The manipulation of the beta functions and phase advances removes the time dependence of the nonlinear kick so as to avoid introducing parametric resonances. In~\cite{Webb2020,Webb2015}, the horizontal and vertical chromaticity should also be equal to maintain integrability for off-momentum particles to their lowest-order approximation.

%\jlnote[inline]{When do we want to spell out "beta" and when do we want to use $\beta$?}
%\jenote[inline]{I would propose that all instances of ``beta functions'' and ``beta profile'' be spelled out, whereas all references to optimization variables or the value at a specific point use $\beta$ symbol. The idea is to avoid compound words with symbols and words. I decided to leave the hyphens off in keeping with typical accelerator stype}
% SW: I agree with above. I did remove the one $\beta$ in this section

The authors of~\cite{EldredIPAC17} and~\cite{EldredIPAC18} present accelerator lattices that combine Danilov--Nagaitsev integrable design criteria with features specific to the Fermilab RCS application. Those lattices were produced with a combination of Nelder--Mead simplex optimization~\cite{NelderMead} and manual manipulation; they serve as a platform for simulation experiments of nonlinear optics. However, they are bespoke lattices that do not represent a thorough or systematic investigation of the optimization landscape. The lattice given in~\cite{EldredIPAC18} was used as a starting point for the optimization in this paper. 

\begin{figure}[tb]
\centering
\includegraphics[width=0.9\linewidth]{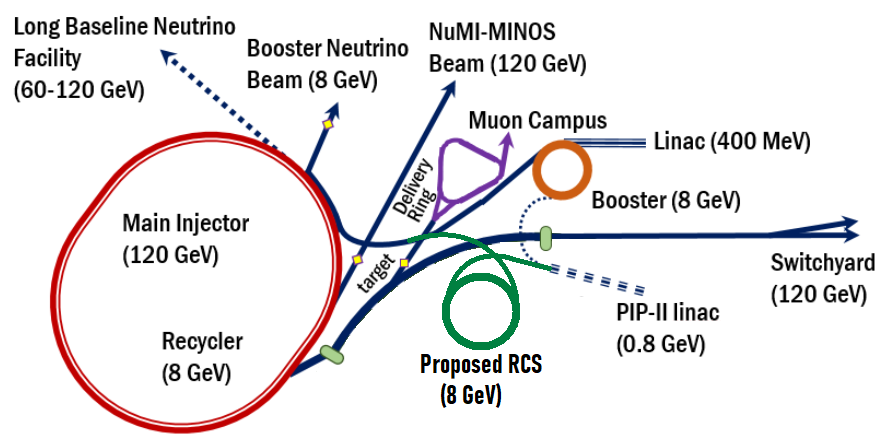}
\caption{PIP-II and LBNF upgrades shown in dashed blue. Proposed RCS upgrade shown in green. Main injector ring (connecting RCS to LBNF) shown in red. Adapted from \cite{Shiltsev2017}.}
\label{fig:RCS_site}
\end{figure}

In addition to the Danilov--Nagaitsev criteria for integrable optics, previous RCS lattice designs were shown to fulfill conventional accelerator optics criteria that are specific to the Fermilab RCS application. The $637.0468718545753$~m circumference 
%GAIL  - I love simple roundoff
of the integrable RCS design was not varied in our optimization; this RCS circumference was previously chosen to satisfy requirements for filling the downstream proton accelerator known as the Fermilab Main Injector. The Main Injector accumulates beam from the RCS, accelerates it to higher energies, and then delivers it to the LBNF beamline; see \Cref{fig:RCS_site}.
%SW: The above is because {fig:RCS_site} is not yet reffed in the text.
The large circumference (relative to the Fermilab Booster) also arises as a result of the requirements for a high number of periodic cells, low-momentum compaction factor, and dispersion-free straights. Increasing the number of periodic cells (around the ring) to twelve was found to improve the performance of the nonlinear integrable optics with intense space-charge by more than a factor of 2. The low-momentum compaction factor ($<5.5 \times 10^{-3}$) is necessary to avoid a loss of longitudinal focusing when accelerating the beam from $\approx$1~GeV to $\approx$8~GeV. The dispersion-free straights eliminate transverse-longitudinal coupling from RF focusing as well as simplify injection and extraction optics.

%\begin{figure}[tbh]
%\centering
%\includegraphics[width=\linewidth]{figures/RCSsite.png}
%\caption{Location of the proposed RCS, relative to the PIP-II linac (which would injects beam into the RCS) and the Main Injector (which would accelerate the beam from the RCS to the LBNF target). \cite{EldredJINST}}
%\label{fig:RCS_site}
%\end{figure}

The RCS dipole magnets are limited to a bending radius of $\rho=21.2$~m, equivalent to a field strength of 1.4~T at 8~GeV, comparable to the Fermilab Main Injector dipoles at 120~GeV. The RCS quadrupole magnets are limited to a focusing strength of $K_{1}=1.2$~m$^{-2}$, a pole radius and pole field comparable to Fermilab Main Injector dipoles at 120~GeV. A minimum of 0.2~m between elements of the magnetic lattice is required to allow clearance for the copper coil winding and vacuum flanges.
These design considerations arise as constraints in our optimization model presented in \Cref{sec:study}.

In~\cite{EldredJINST} it is argued that the maximum values of the beta functions (horizontal and vertical) play a direct role in the cost and performance of an RCS design. The size of the beam is proportional to the square root of the beta functions (neglecting dispersive effects) and the beam emittance (the phase area occupied by the beam). Consequently the beta function determines the minimum diameter of a round beampipe aperture, which is a major consideration in a magnet cost. The physical size of accelerator magnets generally scales quadratically with the diameter of the aperture (transversely, to capture the return flux) or cubically (longitudinally, to maintain the same integrated field). On the other hand, within a given beampipe aperture, a smaller beta function can allow for a proportionately larger beam emittance, which in turn weakens deleterious space-charge effects. Consequently the maximum value taken by the beta functions over the beamline is considered as an objective function for our optimization model.

\begin{figure}[tb]
\centering
\includegraphics[width=0.8\linewidth]{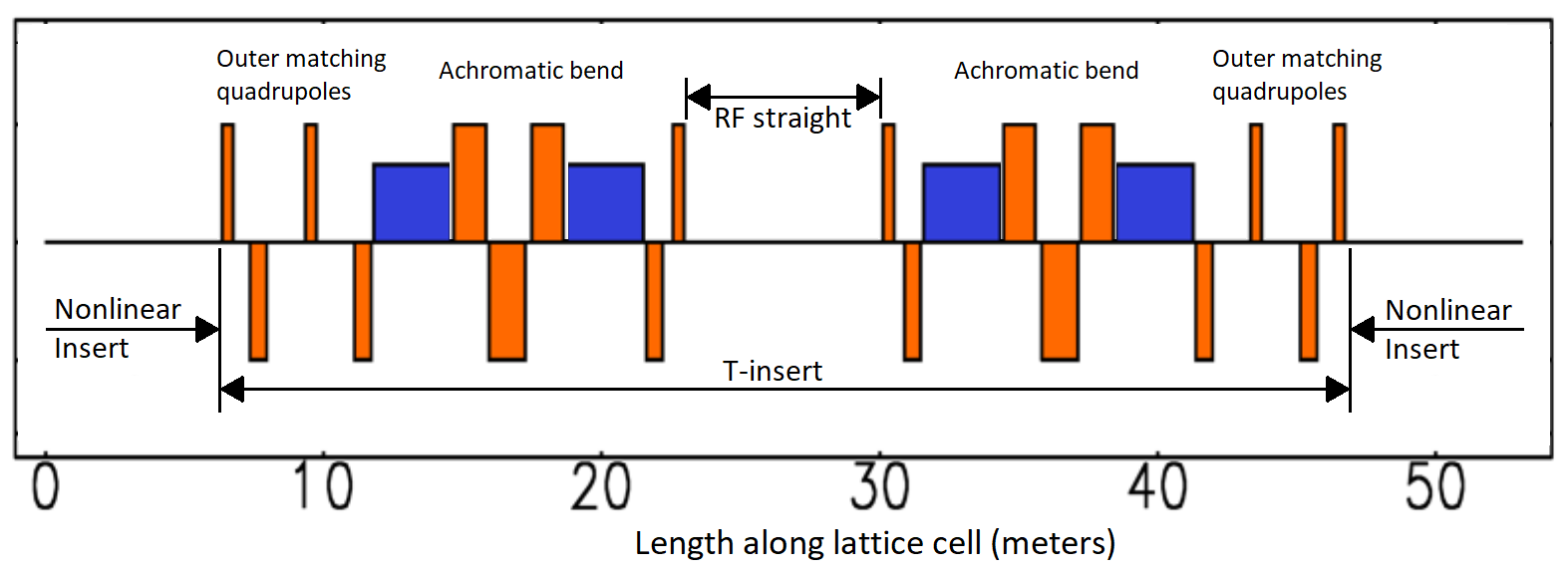}
\caption{Lattice structure for each of 12 RCS cells. Unoptimized location and length of magnetic lattice elements with dipoles shown as short blue rectangles and quadrupoles shown as tall orange rectangles.}
\label{fig:Lattice_Cell}
\end{figure}

As shown in \Cref{fig:Lattice_Cell}, the overall lattice structure of the ring is 12 identical lattice cells, each cell composed of two double-bend achromatic arcs that separate two dispersion-free straight sections of unequal lengths. The straight section in the ``center'' of the cell is for conventional accelerator components---RF accelerating cavities, injection chicane, extraction septa, and collimators. The straight section in the ``outer'' region of the cell is for the nonlinear optics inserts, with all other parts of the cell forming the linear T insert. Each achromatic arc is a quadrupole triplet between two dipole bends and is tightly constrained by the dispersion-matching, low-momentum compaction factor, quadrupole strength, and overall compactness. A simple quadrupole doublet conveys the beam through the ``RF'' straight section. On the outer edges of the cell, at least a quadrupole triplet is needed to convey the beam through the nonlinear straight with matching beta functions; and a fourth quadrupole is included to allow additional finesse between the phase-advance, chromaticity, and maximum beta constraints.

The phase advance over the nonlinear insert is a parameter closely related to the achievable nonlinear tune shift with amplitude as well as the sensitivity of the nonlinear insert to various errors. The IOTA facility, an experimental ring for beam dynamics research~\cite{Antipov2017}, uses a phase advance across the nonlinear insert of $0.3\times2\pi$. The RCS application is focused on robust performance and may allow for a smaller phase advance.

\section{RCS simulation with Synergia}
\label{sec:synergia}

Candidate RCS designs are evaluated by the accelerator modeling
framework Synergia~\cite{Synergia,Amundson:2004qd}, which was developed 
to simulate both extant and proposed accelerator designs, to understand
observed behavior in operational accelerators, to evaluate designs for proposed
upgrades, and to probe beam physics effects~\cite{Macridin:2015vua}.
Synergia combines the two major facets of accelerator physics: single particle
optics and beam collective effects.
To first order, particle transport in a beamline is analogous to light
propagation in an optical system, thus the appellation ``beam optics.''
Collective effects arise due to interactions among all the charged particles in
a beam bunch, primarily space charge, which is mediated by electromagnetic
fields.
Synergia has the capability to track bunches of particles and evaluate
resulting electromagnetic fields using field solvers and applying resulting
forces to the particle motions.
For this study we have neglected collective effects in order
to achieve a first optimization of the RCS design, but the capability in
Synergia leaves open the possibility of considering them in a more detailed
future study.

The particle motion in a periodic focusing system is governed locally by standard equations of motion but can be globally characterized by solutions to 
Hill's equations. 
These are a set of linear
differential equations that depend on the (de)focusing strength as a function
of position.
The focusing strength of an accelerator is determined by beamline element strengths and positions.
These are the decision variables that define the candidate designs described below.
The solutions to Hill's equations consist of beta functions and phase advances (optical functions) that 
appear in the objectives and constraints of our current problem.
%Synergia does not solve Hill's equations explicitly.
%Because the equations are linear, the transformation of particle coordinates can be represented as a linear map, which Synergia derives by propagating the particle coordinates and their first derivatives through the beamline elements using the equations of motion of classical mechanics.
Synergia derives the optical functions from the transfer map calculated by propagating particle coordinates through the beamline elements using the equations of motion of classical mechanics.
%The map of a single pass through the periodic system represented as a matrix is called the \say{one turn map} $M$.
Automatic differentiation of the coordinate vectors produces the maps~\cite{Michelotti2006} that are analyzed to determine stability conditions, beta functions, and phase advances.
There are values for beta functions and phase advances for each beamline
element.
As described in \Cref{sec:rcs}, these quantities make up the bulk of typical  objective and constraint functions.

Synergia is a C++ library for accelerator simulations utilizing
models of physical accelerator components and charged particles organized in
trains of beam bunches along with the machinery to apply the appropriate
physical calculations of particle propagation to the C++ classes and objects. The
computational libraries and data structures are organized so that parallel
computations proceed efficiently on current and future platforms.
The version used for this report uses OpenMPI and OpenMP for parallelism, which
is needed when simulating collective effects among bunches containing
$10^{5}$--$10^7$ particles.
For flexibility and user accessibility, the core objects and methods required
to perform a simulation have Python bindings provided by Boost.Python.
Although simulations may be written in C++, most end users perform Synergia
simulations with scripts written in Python.
Selected classes of user importance may also be extended, including in Python.

This computational architecture proved ideal for this case study.
Candidate RCS designs generated by the optimizer were passed to a Python script, which created the model beamline as Synergia objects using wrapped C++ methods.
Synergia primitives were called to calculate beta functions and other characteristics needed to evaluate the candidate RCS designs.
The results of the evaluation were returned to the optimizer as a Python dictionary.

\section{Mathematical model for the RCS case study}
\label{sec:study}

We now formulate a mathematical optimization problem based on the case study outlined in \Cref{sec:rcs}. 
In particular, we optimize the performance of an integrable lattice by modifying the positions, lengths, and strengths of electromagnetic elements.
Below we detail the key features of the problem: two-sided bound constraints and sequential linear constraints on the elements' positions, lengths, and strengths encoding engineering requirements and nonlinear black-box simulation-based constraints encoding desirable physical properties of the beamline.

\subsection{Decision variables}
\label{sec:vars}

An RCS lattice is defined by the placement of various elements around a ring.
We construct a ring by placing four types of elements (dipoles, quadrupoles,
RF inserts, and NL inserts) on a \textbf{half-cell}, a $\frac{1}{24}$th slice of a full ring. This half-cell is then reflected once (see \Cref{fig:Lattice_Cell}) and repeated twelve times to form the periodic ring. 

We let  
$\x_{e,p,i}$, $\x_{e,l,i}$, and $\x_{e,s,i}$  denote the position, length, and strength, respectively, for the $i$th element of type $e$. We consider $e \in \left\{ d, q, r,n\right\}$, corresponding to  dipole, quadrupole, RF-insert, and NL-insert element types. In each half-cell, we place $\nq = 9$ quadrupoles and $\nd = 2$ dipoles along with a single RF insert and a single NL insert in the fixed ordering
\begin{equation}\label{eq:ordering}
(n,q,q,q,q,d,q,q,q,d,q,q,r),
\end{equation}
as illustrated in \Cref{fig:Lattice_Cell}.

The position (in meters) defines the center of each element. The length (in meters) is the length of the 
element; and thus the starting position of element $e$ is $\x_{e,p,i} - 0.5\x_{e,l,i}$, and its end position is $\x_{e,p,i} + 0.5\x_{e,l,i}$. The quadrupole strength is a scaled focusing strength given in $\text{radians}/\text{m}^2$, and the dipole strength is specified by a bend angle in radians.  
For this optimization, the effects of the RF-insert and NL-insert elements are not simulated, and so they contribute to the lattice design only through the space that they occupy and how the spacing interacts with the constraints on other elements. 

We collect all decision variables in the $(\nx=3\nq+2\nd+1=32)$-dimensional vector $\xb$ with components
\begin{equation}
\label{eq:allx}
\left( 
\x_{q,l,i},
\x_{q,p,i},
\x_{q,s,i}\right)_{i=1}^{\nq},
\left( 
\x_{d,l,i},
\x_{d,p,i}\right)_{i=1}^{\nd},
\x_{d,s,1}.
\end{equation}
Notably, the decision variables we work with do not include five terms ($\x_{n,p}$, $\x_{n,l}$, $\x_{r,p}$, $\x_{r,l}$, and $\x_{d,s,2}$), which are eliminated from our formulation by virtue of additional design considerations discussed next. 
We assume the NL insert starts at the beginning of the half-cell:
\begin{equation}
\label{eq:nnl_starting_pos}
      \x_{n,p} -0.5\x_{n,l} = 0.
\end{equation}
The end of the NL insert coincides with the start of the first quadrupole: 
\begin{equation}
\label{eq:nnl_starting_pos2}
    \x_{q,p,1} - 0.5 \x_{q,l,1} - (\x_{n,p}   + 0.5 \x_{n,l})  = 0.
\end{equation}
There is a fixed circumference $C=$637.0468718545753:
\begin{equation}
\label{eq:Cmax2}
      \x_{r,p}+0.5\x_{r,l} = \frac{C}{24}.
\end{equation}
Since the ring must close, the bend angle in a half-cell must be $\pi/12$:
\begin{equation}
\label{eq:bend_angles}
        \x_{d,s,1}+\x_{d,s,2} =\frac{\pi}{12};
\end{equation}
and the start of the RF insert must coincide with the end of the last quadrupole:
\begin{equation}
    \x_{r,p} - 0.5 \x_{r,l} - (\x_{q,p,9} + 0.5 \x_{q,l,9})  = 0.
    \label{eq:rf_quad_gap}
\end{equation}
These linear equalities allow us to eliminate the
terms $\x_{n,p}$, $\x_{n,l}$, $\x_{r,p}$, $\x_{r,l}$, and $\x_{d,s,2}$, which can be recovered via
\begin{align*}
    \x_{n,p} &= 0.5(\x_{q,p,1} - 0.5\x_{q,l,1})
    \\
    \x_{n,l} &= \x_{q,p,1} - 0.5\x_{q,l,1}
    \\
    \x_{r,p} &= \dfrac{C}{24} - 0.5\left(\dfrac{C}{24} - \x_{q,p,9} - 0.5\x_{q,l,9}\right)  
    \\
    \x_{r,l} &= \dfrac{C}{24} - \x_{q,p,9} - 0.5\x_{q,l,9}
    \\
    \x_{d,s,2} &= \dfrac{\pi}{12}  - \x_{d,s,1}.
\end{align*}

Following the taxonomy in~\cite{taxonomy15}, we separate constraints on the remaining $\xb$ into simulation-based constraints (i.e., those for which
(in)feasibility can be verified without a Synergia call) and algebraic constraints (i.e., those that do not require a Synergia call).

\subsection{Algebraic constraints}
\label{sec:algcon}

Our model's algebraic constraints consist of two-sided bound constraints and 
linear inequality constraints. Together, these algebraic constraints define a compact polytope in $\R^{\nx}$.

Two-sided bound constraints (with unequal lower and upper bounds)
have the benefit of not contributing twice to the combinatorial complexity of an optimization problem:
at most one of the sides can be active, and thus the two bounds will not contribute to potentially violating a constraint qualification  (e.g., LICQ) unless the objective's derivative with respect to a decision variable vanishes.
We assume that the $\nq+\nd$ lengths are bounded:
\begin{align*}
      0 \leq \x_{q,l,i} &\leq 5, \qquad i=1,\ldots,\nq\\
      \x_{d,l,i} &\leq 5, \qquad i=1,\ldots,\nd.
\end{align*}
We also assume that the $\nq+1$ strengths are bounded. By noting that the lower bound on the dipoles can be
tightened, since the dipole strengths must obey the bend angle constraint in
\cref{eq:bend_angles}, we arrive at
\begin{align*}
      -1.2 &\leq \x_{q,s,i} \leq 1.2, \qquad i=1,\ldots,\nq \\
      \pi/12 - 0.2    &\leq \x_{d,s,1} \leq 0.2.
\end{align*}

The linear inequality constraints encode the fixed relative positioning of the elements in the order given in \cref{eq:ordering} with spacing, as well as the reduced circumference and bend angle constraints. We have
\begin{align*}
    \x_{q,p,1}   - 0.5\x_{q,l,1} &\ge 1 \\
    \x_{q,p,i} - 0.5 \x_{q,l,i} - (\x_{q,p,i-1} + 0.5 \x_{q,l,i-1}) &\ge 0.2, \qquad i=2,\ldots,4 \\
    \x_{d,p,1} - 0.5 \x_{d,l,1} - (\x_{q,p,4} + 0.5 \x_{q,l,4}) &\ge 0.2 \\
    \x_{q,p,5} - 0.5 \x_{q,l,5} - (\x_{d,p,1} + 0.5 \x_{d,l,1}) &\ge 0.2 \\
    \x_{q,p,i} - 0.5 \x_{q,l,i} - (\x_{q,p,i-1} + 0.5 \x_{q,l,i-1}) &\ge 0.2, \qquad i=6,7 \\
    \x_{d,p,2} - 0.5 \x_{d,l,2} - (\x_{q,p,7} + 0.5 \x_{q,l,7}) &\ge 0.2 \\
    \x_{q,p,8} - 0.5 \x_{q,l,8} - (\x_{d,p,2} + 0.5 \x_{d,l,2}) &\ge 0.2 \\
    \x_{q,p,9} - 0.5 \x_{q,l,9} - (\x_{q,p,8} + 0.5 \x_{q,l,8}) &\ge 0.2.
\end{align*}
We note that the first constraint combines \cref{eq:nnl_starting_pos} with \cref{eq:nnl_starting_pos2} and a lower bound of 1 on the NL-insert length. 

We include two additional linear inequalities that combine the bounds $[3,3.75]$ on the RF-insert lengths, the fixed circumference constraint \cref{eq:Cmax2}, and the spacing constraint \cref{eq:rf_quad_gap} between the last quadrupole and RF insert:
\begin{align}
   \x_{q,p,9} + 0.5 \x_{q,l,9} &\le \dfrac{C}{24} - 3
   \\
    -\x_{q,p,9} - 0.5 \x_{q,l,9} &\le 3.75 - \dfrac{C}{24}.
\end{align}
Lastly, we apply two linear constraints to proportionally lower bound the dipole lengths, $\rho\x_{d,s,i} \leq \x_{d,l,i}$ for $i=1,\ldots,\nd$. By using the change of variables defined in \cref{eq:bend_angles}, we have
\begin{align}
    \rho\x_{d,s,1} - x_{d,l,1} &\le 0 \\
    - \rho\x_{d,s,1} - x_{d,l,2} &\le -\dfrac{\rho\pi}{12},
\end{align}
where $\rho = 21.2$ meters. 

These algebraic constraints (i.e., the $2\nq+1$ two-sided bound constraints, $n_d$ one-sided bound constraints, and the 15 other linear inequalities) can be expressed compactly by $\Ab \xb \leq \bb$, where $\Ab$ is of size 
 $\nac = 4\nq+19$ by $\nx=3\nq+2\nd+1$ (55 by 32), which we use for conciseness below.
 
For the remainder of this study we rescale the decision variables $\xb$ to the unit cube. The appropriate rescaling is by the side lengths of the smallest box that encloses the polytope $\left\{ \xb \in \R^{\nx} : \, \Ab \xb \leq \bb \right\}$.

\subsection{Simulation-based constraints}
\label{sec:simcon}

Numerous criteria (described in \Cref{sec:rcs}) are desired in an RCS lattice. 
In our model we opt for a constraint-based approach in which simulation-based constraints encode the majority of the physical design criteria for the RCS lattice described in \Cref{sec:rcs}. With this approach a key concern is establishing joint feasibility across all the constraints;
we are also interested in the dependence of solutions on the constraint parameters, which we study in \Cref{sec:kkt}. 
In the remainder of this section we describe these $\nsc = 23$ constraints, all of which require a Synergia call in order to verify (in)feasibility. 

A fundamental constraint to the RCS design and mathematical model is the requirement that the lattice have a periodic orbit (and is therefore stable). Many quantities of interest are undefined for physically infeasible solutions where no periodic orbit exists.  In those situations Synergia is unable to evaluate the candidate solution.  Mathematically, the periodic orbit exists when the \say{one turn map} $M(\xb)$, corresponding to direction of motion ($x$, $y$, $z$), has eigenvalues 
of roughly unit norm. 
We constrain only the transverse directions $x$ and $y$ so we require two conforming eigenvalue pairs.  We enforce the norm unity conditions with the two two-sided constraints $c_1,c_{2}$ on the modulus of the single relevant eigenvalue $\lambda(\xb)$ of $M(\xb)$: 
\begin{align}
     c_1(\xb) &:= | \lambda(\xb) | - 1-\epsilon_\lambda \leq 0
     \\
     c_2(\xb) &:= -| \lambda(\xb) | + 1-\epsilon_\lambda \leq 0,
\end{align}
with $\epsilon_\lambda = 10^{-10}$. 
In practice, we find that this constraint is 
always satisfied at $\xb$ that are feasible with respect to all other constraints. 

The next constraint enforces that the momentum compaction factor $\alpha_{c}(\xb)$ should be at most $5.5 \times 10^{-3}$:
\begin{align}
c_{3}(\xb)  := \alpha_{c}(\xb) - 5.5\times10^{-3} \leq 0.
\end{align}
The RCS design seeks a dispersion-free RF insert, which in practice is relaxed by a numerical tolerance $10^{-3}$:
\begin{align}
   c_{4}(\xb)  &:= D_{x,\text{rf}}(\xb) - 10^{-3}  \leq 0\\
   c_{5}(\xb)  &:= -D_{x,\text{rf}}(\xb) -10^{-3} \leq 0.
\end{align}
Similarly, constraints for dispersion-free NL inserts are relaxed by a numerical tolerance $10^{-5}$:
%\Psi_{x,l} -> D_{x,\text{nl1}}
%\Psi_{x,l} -> D_{x,\text{nl2}}
\begin{align}
    c_{6}(\xb)  &:=D_{x,\text{nl1}}(\xb) - 10^{-5} \leq 0\\
    c_{7}(\xb)  &:=-D_{x,\text{nl1}}(\xb) - 10^{-5} \leq 0\\
    c_{8}(\xb)  &:=D_{x,\text{nl2}}(\xb) - 10^{-5} \leq 0\\
    c_{9}(\xb)  &:=-D_{x,\text{nl2}}(\xb) - 10^{-5} \leq 0.
\end{align}
To ensure sufficient similarity of the beta functions over the NL insert, we enforce the four constraints
\begin{align}
c_{10}(\xb) &:= \beta_{y,1}(\xb) - (1+\epsilon_\beta)\beta_{x,1}(\xb)\leq 0 
\\
c_{11}(\xb) &:= -\beta_{y,1}(\xb) + (1-\epsilon_\beta)\beta_{x,1}(\xb)\leq 0  
\\
c_{12}(\xb) &:= \beta_{y,2}(\xb) - (1+\epsilon_\beta)\beta_{x,2}(\xb)\leq 0
\\
c_{13}(\xb) &:= -\beta_{y,2}(\xb) + (1-\epsilon_\beta)\beta_{x,2}(\xb)\leq 0, 
\end{align}
where we use $\epsilon_\beta = 0.01$ and where $\beta_{\cdot,1}$ indicates the beta function value at the start of the NL insert and $\beta_{\cdot,2}$ indicates its value at the end of the NL insert. 

We enforce integrability constraints by using
the modulo operator to ensure that the T-insert phase advances $\Phi_x$ and $\Phi_y$ are multiples of $\pi$. We apply a shift of $\pi/2$ to move away from the discontinuity in the $\bmod(\pi)$ function:
\begin{align}
  c_{14}(\xb) &:= (\pi/2 + \Phi_x(\xb))\bmod(\pi)  - (\pi/2+\epsilon_{\Phi}) \leq 0
  \\
  c_{15}(\xb) &:= -(\pi/2 + \Phi_x(\xb))\bmod(\pi) + (\pi/2-\epsilon_{\Phi} ) \leq 0
  \\
  c_{16}(\xb) &:= (\pi/2 + \Phi_y(\xb))\bmod(\pi)  - (\pi/2+\epsilon_{\Phi}) \leq 0
  \\
  c_{17}(\xb) &:= -(\pi/2 + \Phi_y(\xb))\bmod(\pi) + (\pi/2-\epsilon_{\Phi}) \leq 0,
\end{align}
where we use the tolerance $\epsilon_{\Phi}=10^{-4}$.

We also seek an RCS design where the ring-wide betatron tunes $\nu_x(\xb),\nu_y(\xb)$ are sufficiently close. Here we use a tolerance of $10^{-4}$ and thus have
\begin{align}
  c_{18}(\xb) &:= \nu_x(\xb) - \nu_y(\xb) - 10^{-4} \leq 0
  \\
  c_{19}(\xb) &:= - (\nu_x(\xb) - \nu_y(\xb)) - 10^{-4} \leq 0.
\end{align}
The phase advance over the NL inserts $\Psi_{x}$ and $\Psi_{y}$ should be at least $0.2\times 2\pi$:
\begin{align}
  c_{20}(\xb) &:=-\Psi_x(\xb) + 0.2\times2\pi \leq 0 \\
  c_{21}(\xb) &:=-\Psi_y(\xb) + 0.2\times2\pi \leq 0.
\end{align}
The horizontal and vertical chromaticity 
should differ by no more than $0.1$ globally:
\begin{align}
c_{22}(\xb) &:=C_{x}(\xb)-C_{y}(\xb) - 0.1 \leq 0
\\
c_{23}(\xb) &:= -(C_{x}(\xb)-C_{y}(\xb)) - 0.1 \leq 0.
\end{align}

Collectively, we denote these 23 simulation-based constraints by $\cb(\xb) \leq \zerob$, so that nonpositive values indicate feasibility. 

With a complete description of the decision variables and constraints in hand, we have a mathematical description of our feasible decision space, which aims to capture RCS lattice designs satisfying the intent of \Cref{sec:rcs}. In the following section we summarize challenges for finding feasible decision points and establish that such feasible decisions exist.

\section{Phase-one solution: finding a feasible point}
\label{sec:feasibility}

Most optimization approaches to problems of the form \cref{eq:genprob} benefit from being initialized with a feasible starting point. Nonlinear simulation-based constraints can make finding such a point difficult, and 
\say{phase-one} (i.e., feasibility seeking) optimization approaches may not cope well with the nonlinear, black-box form of the constraints. Finding a feasible point is a key challenge of this problem since it validates (up to simulation fidelity) the existence of an RCS lattice with the desired properties outlined in \Cref{sec:rcs}, which has until now not been shown to exist. 
In our search for a feasible point we take advantage of the polytope structure in our algebraic constraints (since this is available without querying the simulation) and a solution that is known to satisfy some (but not all) of the constraints from previous synchrotron design problems.

The polytope defined by $\Ab\xb \leq \bb$ bounds the feasible region $\Omega$, and thus sampling within the polytope may appear to be a promising method of finding feasible points. Approaches for random sampling from the polytope include taking convex combinations of the polytope vertices or employing Markov chain Monte Carlo (MCMC) methods such as hit-and-run sampling or the Vaidya walk~\cite{chen2018polytope,mete2012pattern}. MCMC methods require only an interior point to the polytope. Although vertex enumeration is difficult in general, it is tractable in our case because our linearly constrained polytope has some dimensions (e.g., the element strengths) that appear only in bound constraints. This allows an easier enumeration of vertices, for example, by enumerating the vertices of the bound-constrained region separately from the other vertices of the linearly constrained region~\cite{bueler2000vinci,Khachiyan2008}. 

\begin{figure}[tb]
\centering
\includegraphics[width=0.8\linewidth]{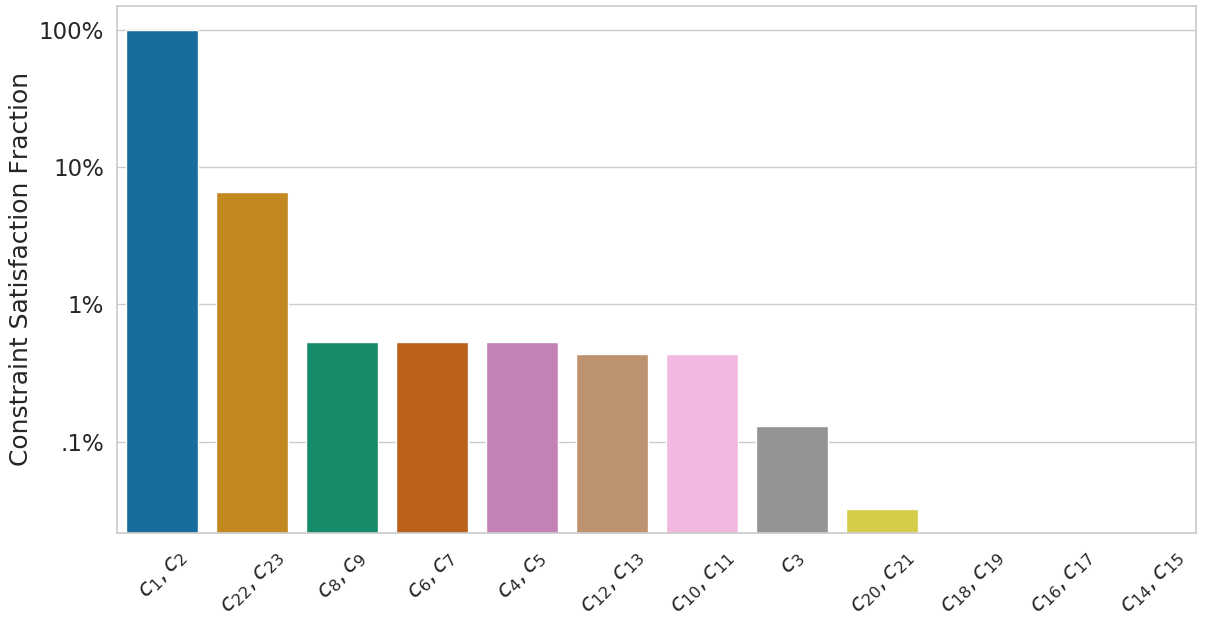}
\caption{Fraction of $6,182$ successfully simulated points randomly sampled within the polytope at which each simulation-based constraint (or both sides of a two-sided constraint pair) is satisfied. The two-sided constraint pairs $(c_{18},c_{19})$, $(c_{16},c_{17})$, and $(c_{14},c_{15})$ are never satisfied.}
\label{fig:feasibility_fraction}
\end{figure}

\Cref{fig:feasibility_fraction} illustrates 
that taking uniformly selected convex combinations of vertices is not a useful approach in finding feasible points in our case study: in 67,810 samples, not a single point was feasible. Of the 67,810 points sampled, roughly 90\% resulted in unsuccessful simulation evaluations (i.e., because the stability conditions $c_1,c_2$ were not satisfied). \Cref{fig:feasibility_fraction} shows that even the successfully evaluated points still tend to violate most simulation-based constraints. Furthermore, MCMC methods were observed to have an even higher fraction of unsuccessful simulation evaluations than sampling convex combinations of vertices have. 
We conclude that the points that are feasible for the simulation-based constraints are \say{needles in the haystack} of the points satisfying the linear constraints.

\begin{figure}[tbh]
    \centering
    \includegraphics[width=.8\linewidth]{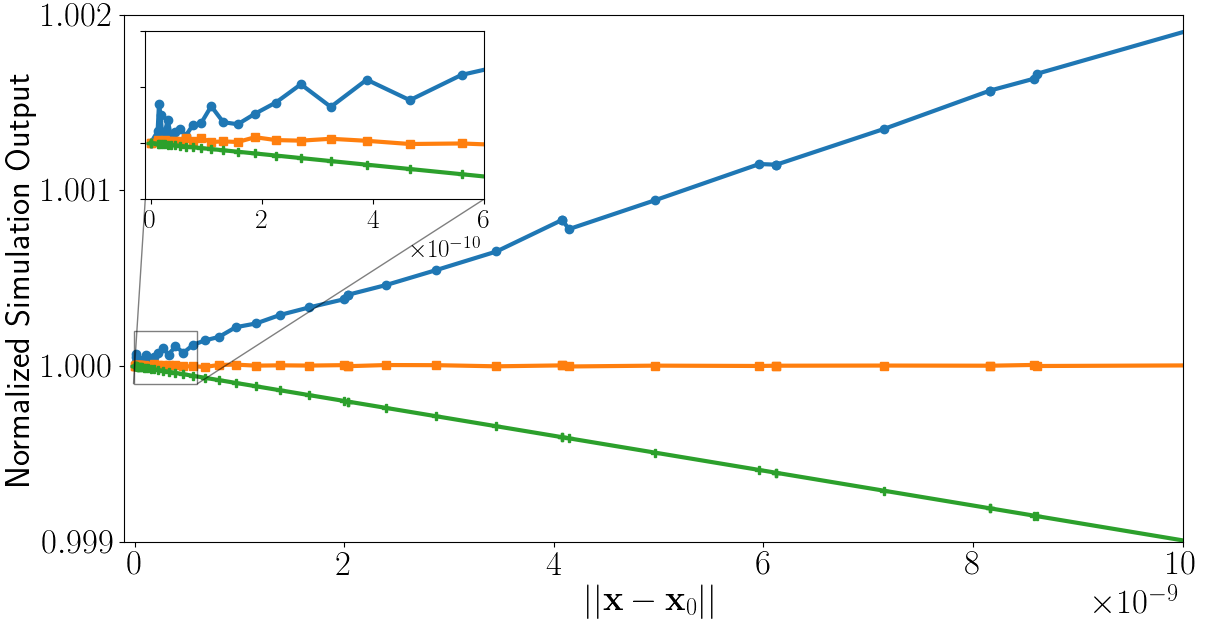}
    \caption{Examples of smoothness in simulation outputs momentum compaction $\alpha_c$ (blue), NL-Dispersion $D_{x,\text{nl1}}$ (orange) and betatron tunes $\nu_x - \nu_y$ (green) for changes in $\xb$ smaller than $10^{-9}$ along a ray in parameter space. The decision vectors $\xb$ and $\xb_0$ are on the unit cube. Notably, Synergia shows high fidelity under small changes in the decision variables $\xb$, free of any significant noise to roughly a $10^{-10}$ change in $\xb$. In the vicinity of feasible solutions there are no nonlinearities in the optics that would lead to chaotic behavior.}
    \label{fig:noise_in_sim}
\end{figure}

On the other hand, a phase-one optimization procedure initialized at an (infeasible) \say{expert-designed} point, denoted $\eb_0$, was successful in finding a feasible point. The point $\eb_0$ was lacking for two reasons. First, it violated algebraic constraints: the spacings between Dipole 2 and Quadrupole 7, Dipole 2 and Quadrupole 8, Dipole 1 and Quadrupole 4, Dipole 1 and Quadrupole 5, and Quadrupoles 5 and Quadrupole 6 had a separation distance of 0.1~m rather than the desired 0.2~m. Second, the simulation-based constraints $c_7,c_9,c_{14},c_{16},c_{19}$ were violated; see \Cref{fig:lollipop_initial_optima}. To ameliorate these issues, we solved two optimization problems: one to move into the algebraic-constraint-defined polytope while not worsening violation of the simulation-based constraints and a subsequent one to improve the violation of the simulation-based constraints. 
Working with simulation-based constraints can be challenging because of discontinuities or numerical noise in a constraint function, especially when these are present when a constraint is nearly active  (see, e.g., \cite[Section 7]{LMW2019AN}, \cite[Chapter 12]{AudetHare2017}).
In our problem, we found that the simulation outputs were remarkably smooth; for example, as shown in \Cref{fig:noise_in_sim}, we see that numerical noise is not readily apparent for changes in $\xb$ larger than $10^{-10}$.
This is expected since the system being simulated is a smooth physical system with forces that vary only linearly with particle amplitudes, which should not exhibit chaotic behavior with small deviations.

Starting from $\eb_0$ and letting $\Ab_i$ indicate the $i$th row of $\Ab$,
we minimize the algebraic constraint violation $\sum_{i=1}^{\nac}\max\{\Ab_i\xb - b_i,0\}^2$ while constraining the simulation-based constraints and algebraic constraints that hold at $\eb_0$ to continue to hold and constraining the simulation-based constraints that are violated at $\eb_0$ by the initial value $c_i(\eb_0)$ to prevent them from worsening. A globally optimal solution to this problem (i.e., satisfying the stated constraints and all algebraic constraints) was found by using the COBYLA local optimization routine~\cite{Powell1994,pdfo}. 

From this point, $\eb_1$, in the polytope we initialize a second optimization problem in order to find simulation-based constraint satisfaction. We minimize the simulation-based constraint violation $\sum_{i=1}^{\nsc}\max\{c_i(\xb),0\}^2$ while constraining the simulation-based constraints and algebraic constraints that hold at $\eb_1$ to continue to hold and constraining the simulation-based constraints $c_i$ that are violated at $\eb_1$ by their initial value $c_i(\eb_1)$ to prevent their infeasibility from worsening. We also constrain the values of all components of $\betab_x$ and $\betab_y$ by the nominal value of $50$ so as to not unnecessarily affect objective function values (discussed in the next section) when searching for a feasible point. COBYLA solved this problem to global optimality (i.e., finding a point in $\Omega$ that also satisfied constraints imposed on $\betab_x$ and $\betab_y$) in roughly 60,000 function evaluations.

A key to the success of the phase-one optimization was starting the feasibility restoration at an expert-designed point that was relatively close to the feasible region $\Omega$. In part because of the nonconvexity of the simulation-based constraints, we found that beginning phase-one from points with significant simulation-based constraint violation did not produce feasible solutions.

\section{Nondifferentiable composite objective functions and local optimality}
\label{sec:objective}

Finding points satisfying all of the constraints posed in \Cref{sec:study} allows us to consider minimization of the beta functions. A minimal beta profile is desired because it indicates a high level of beam concentration.

Because the accelerator community does not have a single standard objective function, we investigate the solution of \cref{eq:genprob} with three different objectives, consisting of different nonsmooth compositions of the simulation outputs $\betab_x(\xb)$ and $\betab_y(\xb)$. 
Expert solicitation indicated that a leading criterion is a low maximum in the beta functions. A second criterion is a symmetry in the peaks of the beta functions. 
Based on these criteria and a composite objective of the form $f(\xb) = h(\betab(\xb))$, we consider the following three compositions:
\begin{eqnarray}
   h_a(\betab(\xb)) &=&  \max_j \{\beta_{x,j}(\xb)\}
   + \max_j \{\beta_{y,j}(\xb)\} = \left\| \betab_x(\xb) \right\|_{\infty} + \left\| \betab_y(\xb) \right\|_{\infty} \qquad
   \label{eq:fa}
\\
   h_b(\betab(\xb)) &=&  \max \left\{ \max_j \{\beta_{x,j}(\xb)\}
   , \max_j \{\beta_{y,j}(\xb)\} \right\} = \left\| \betab(\xb) \right\|_{\infty}
   \label{eq:fb}
\\
h_c(\betab(\xb)) &=&\max_{j}\left\{ \beta_{x,j}(\xb)+\beta_{y,j}(\xb)\right\} = \left\| \betab_x(\xb) +  \betab_y(\xb) \right\|_{\infty}.
   \label{eq:fc}
\end{eqnarray}
Since the beta functions are nonnegative, we have that,
for any $\xb \in \Omega$, $h_b(\betab(\xb)) \leq h_c(\betab(\xb)) \leq h_a(\betab(\xb))$. 
Although smooth approximations to these objectives can be formulated, we directly employ these nonsmooth forms in order to preserve the interpretability of each objective form.

\begin{figure}[tb]
    \centering
    \includegraphics[width=.6\linewidth]{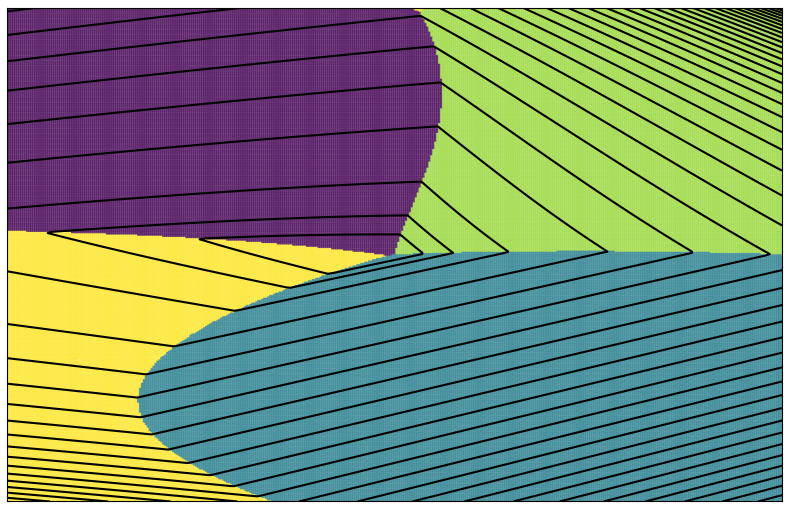}
    \caption{Visible nonsmoothness of the objective defined by $h_b$ in a random two-dimensional slice of $\R^{\nx}$ around an arbitrary feasible point. In this plot, $h_b$ is defined by the active indexes $\beta_{x,48}$ (bottom left), $\beta_{x,180}$ (top right), and $\beta_{y,111}$ (bottom right), $\beta_{y,158}$ (top left).
    }
    \label{fig:nonsmooth_objectives}
\end{figure}

\Cref{fig:nonsmooth_objectives} illustrates that the resulting objective functions readily exhibit nondifferentiable behavior in $\Omega$.

\subsection{Manifold sampling from feasible points}
\label{sec:ms}

To address the structured nondifferentiability illustrated in \Cref{fig:nonsmooth_objectives}, we employ a manifold sampling trust-region method~\cite{KLW18,LMW16,Larson2020}.
Since the nonsmoothness occurs because of a composition of a known mapping $h$ with a (smooth) simulation output $\betab(\xb)$, manifold sampling can exploit the known form of each $h$. Manifold sampling is a model-based approach that categorizes points $\betab(\xb)$ in the domain of $h$ as belonging to different manifolds that occur due to the nonsmoothness in $h$; this information is used when determining search directions~\cite{KLW18,LMW16,Larson2020}.

Instead of enforcing all 78 inequality constraints, we  use a quadratic penalty formulation~\cite{nocedal2006no}, where we penalize the square of the constraint violation $[\cdot]_+ = \max\{\cdot,0\}$ through minimization of the unconstrained objective 
\begin{equation}
f(\xb) = h(\betab(\xb)) + \alpha \sum_{i=1}^{\nsc} [c_i(\xb)]_+^2 + \alpha \sum_{i=1}^{\nac}\left[\Ab_i\xb-b_i\right]_+^2
\label{eq:penalizedobj1}
\end{equation}
for the various nonsmooth definitions of $h$ above. A sequence of solutions of \cref{eq:penalizedobj1} as $\alpha\to\infty$ converge to a  Karush–-Kuhn–-Tucker (KKT) point of \cref{eq:genprob}, where $f(\xb) = h(\betab(\xb))$ provided that certain conditions are satisfied: feasibility of the limit point, appropriate constraint qualification, and sufficient accuracy of solutions of \cref{eq:penalizedobj1}~\cite[Thm.~17.2]{nocedal2006no}.

Manifold sampling can address the fact that the $f$ in \cref{eq:penalizedobj1} is a nonsmooth function with additional potential nonsmoothness coming from constraint violations. For example, if the constraint and beta functions are continuously differentiable, the Clarke subdifferential of \cref{eq:penalizedobj1} when $h = h_b$ can be computed by taking the convex hull of the gradients of the active beta functions 
\begin{equation}
\begin{aligned}
\partial h_b(\betab(\xb)) =& \underset{j\in J^*(\xb)}{\mbox{co}}   \left\{\nabla\beta_{j}(\xb) \right\} \\
&+ 2 \alpha \left(\sum_{i=1}^{\nsc}[c_i(\xb)]_+\nabla c_i(\xb) + \sum_{i=1}^{\nac}\left[\Ab_i\xb-b_i\right]_+ \Ab_i^T\right),
\end{aligned}
\end{equation}
where $J^*(\xb)= \arg\max_{j}\left\{ \beta_{j}(\xb) \right\}$. It is precisely the knowledge of which manifolds (in this case, the indices of the beta functions) are active that allows manifold sampling to approximate this subdifferential and determine descent directions. Since the quantities $\nabla_{\xb}\beta_{x,\cdot}$ and $\nabla_{\xb}\beta_{y,\cdot}$ are unavailable, 
manifold sampling builds local models $m_{\beta_j}$ of both beta functions at each of a finite number of indices $j$ and uses $\nabla_{\xb} m_{\beta_j}$ in place of the corresponding (unavailable) gradients.

Manifold sampling notes what indices define the set $J^*(\xb)$ at (or near) candidate 
points $\xb$ and uses this information in order to infer when the subdifferential of $f$
changes. Given a set $\Yb$ of past points evaluated during the course of the algorithm (either iterates of the algorithm or from building models of the components of $\betab$) and a trust-region radius $\Delta$, 
the manifold sampling implementation used in our numerical tests approximated the set $J^*(\xb)$ by 
$$\displaystyle \bigcup_{\yb \in \Yb : \| \yb - \xb \| \le \Delta } \left\{ j \, : \, \beta_{j}(\yb) \ge \max_i \left\{\beta_{i}(\yb)\right\} - 10^{-8} \ \right\}.$$
That is, the indices that are considered to define $h_b(\betab(\yb)$ are those where the max is within $10^{-8}$ of the true maximum; we take the union of these sets for any point $\yb$ in the trust region around $\xb$. 

Such an approach to approximating the subdifferential at $\xb$ may seem excessive,
but it is in fact essential for algorithmic performance. The logic of the manifold sampling algorithm requires that the manifolds active at any putative iterate be known before a step can be taken/rejected. This requires an inner set of 
manifold sampling iterations where the set of active manifolds grows iteratively. This process is guaranteed to terminate because the number of indices that can be active is finite. Warm-starting the manifold sampling process with past information means the algorithm can make progress sooner.

\subsection{Local optimality conditions}
\label{sec:kkt}
A key benefit of optimization models in the RCS design is providing insight into how the performance of the RCS design is dependent on tunable parameters. From a physical standpoint this brings interpretability to how design adjustments alter performance.
We are particularly interested in the effect of parameters $\epsb\in\R^{\nsc}$ defining a feasible region based on a parameterized right-hand side of the simulation-based constraints:
\begin{equation}
\Omega_{\epsb} = \left\{\xb\in \R^{\nx} : \Ab \xb \leq \bb, \; c_i(\xb) \leq \eps_i, \, i =1,\ldots, \nsc\right\};    
\label{eq:epsfeasible}
\end{equation}
$\eps_i=0$ corresponds to our nominal formulation, whereas a positive (negative) $\eps_i$ value corresponds to a relaxation (tightening) of the $i$th simulation-based constraint.
We look at the KKT optimality conditions~\cite{nocedal2006no}, which indicate how changes in the violation of active constraints locally improve or degrade performance in the objective. Solving the KKT conditions yields Lagrange multipliers $\lambda_i$ for each constraint; these multipliers quantify how each parameter $\eps_i$ affects the optimal objective value. With the Lagrange multipliers in hand, an RCS designer can weigh tradeoffs between constraint violation and objective quality when considering different designs.

Assuming that the constraint and beta functions are continuously differentiable, the nonsmooth KKT conditions~\cite{bagirov2014intrononsmooth,Fletcher} state that under a constraint qualification, such as LICQ, there exist a subgradient $ \subg(\xb) \in \partial h(\betab(\xb))$ and Lagrange multipliers $\lambdab = (\lambdab^s,\lambdab^a)$ such that
\begin{align}
& \subg(\xb) +\sum_{i\in\mathcal{A}_{\rm sc}(\xb)}\lambda_i^s\nabla c_i(\xb) +\sum_{j\in\mathcal{A}_{\rm ac}(\xb)}\lambda^a_j \Ab_j^T = 0
\label{eq:kkt_stationarity}
\\
&c_i(\xb) \leq \eps_i, \ \ \ i=1,\ldots,\nsc
\\
&\Ab\xb \leq \bb
\\
&\lambdab^s,\lambdab^a \geq \zerob
\\
&\lambda_i^s c_i(\xb) = 0, \ \ \  i = 1,\ldots,\nsc
\\
&\lambda_i^a \left(\Ab_i\xb-b_i\right) = 0, \ \ \  i = 1,\ldots,\nac,
\end{align}
where $\mathcal{A}_{\rm ac}(\xb)$ and $\mathcal{A}_{\rm sc}(\xb)$ indicate the active sets of algebraic and simulation-based (i.e., $c_i(\xb) = \eps_i$) constraints, respectively. 
The stationarity condition \cref{eq:kkt_stationarity} provides us the handle for understanding parametric dependence since it expresses the objective gradient in terms of active constraint gradients. We can use this to approximate the objective value under changes to constraint values, up to linearization. For example, suppose there exists (in fact there does exist by LICQ) a direction $\db$ that is orthogonal to all active constraint gradients aside from, say, $\nabla c_i(\xb)$. Then, the change in objective value in the direction $\db$ is approximately $- \lambda_i^s\db^T\nabla c_i(\xb)$.

When considering the objective $h = h_b$, the subgradients are defined by the convex combination $\subg(\xb) =\sum_{j\in J^*(\xb)} \gamma_j\nabla\beta_{j}(\xb)$, where $\gammab\geq \zerob$ satisfies $\sum_j \gamma_j = 1$
and $J^*(\xb)= \arg\max_{j}\left\{ \beta_{j}(\xb) \right\}$ indexes the set of active beta indices. When $\xb$ is optimal for \cref{eq:genprob}, the Lagrange multipliers $\lambdab$ solve 
\begin{align}
\label{eq:solve_for_lagrange}
\min_{\lambdab, \gammab}\ \ \  \left \|\sum_{j\in J^*(\xb)} \gamma_j\nabla \beta_{j}(\xb) \right. &  \left. + \sum_{i\in\mathcal{A}_{\rm sc}(\xb)}\lambda^s_i\nabla c_i(\xb) +\sum_{i\in \mathcal{A}_{\rm ac}(\xb)}\lambda^a_i \Ab_i^T \right \|^2
\\
\lambdab,\gammab &\geq \zerob 
\\
\lambda^s_i (c_i(\xb)-\eps_i) &= 0  \quad i=1,\ldots,\nsc
\\
\lambda^a_i \left(\Ab_i\xb - b_i\right) &= 0  \quad i=1,\ldots,\nac
\\
\sum_{j\in J^*(\xb)} \gamma_j &= 1.
\end{align}
We use this method to compute $\lambdab$ in \Cref{sec:results-tradeoffs}. In practice, we determine the active set of beta indices to be those functions $\beta_j$ within $10^{-8}$ of the max objective value. The active set of constraints is determined by first finding the steepest descent direction and then collecting the constraints that are violated when traveling in such a direction. In this way we find a sequential list of constraints that are limiting motion. The constraint right-hand sides are then adjusted so that the active constraints exactly equal zero at the optima.
%SW: I bet the last sentence needs clarification

\section{Numerical results}
\label{sec:results}
Here we numerically explore properties of the optimization model \cref{eq:genprob}, described in \Cref{sec:rcs} and developed in \Cref{sec:study}. We discuss the physical properties, physical validity, sensitivity, and effects of integrability at various local optima and compare them with unoptimized points. Furthermore, to understand the dependence of optima on constraint parameters, we analyze the KKT conditions and Lagrange multipliers, which provide a lens to elucidate the dependence of the objective value on the set of restrictive constraints. The principal finding is that the thorough development of an optimization model provides a powerful tool for finding high-performance RCS lattice designs that meet the desired criteria as well as analyzing the design criteria and their effect on solutions. 

\subsection{Studying different objectives}
\label{sec:results-obj}
\begin{figure}[tb]
\centering
\includegraphics[width=0.8\linewidth]{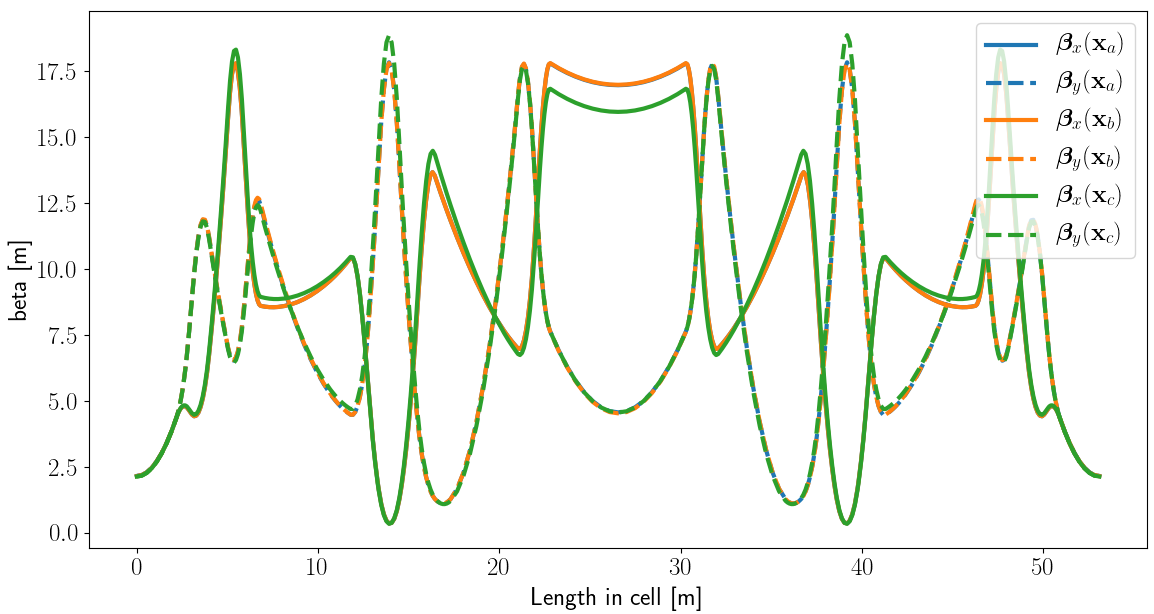}
\caption{Beta functions at the local minima $\xb_a,\xb_b,\xb_c$ after manifold sampling optimization under the respective objectives $h_a,h_b,h_c$. The $h_b$ and $h_c$ objectives produce minima with nearly identical beta functions; the $h_a$ minimizer also has a similar profile.}
\label{fig:three_beta_plots}
\end{figure}

Since it is not immediately clear which beta profiles are preferable, we optimize all three composite objective functions $h_a(\betab(\xb)),h_b(\betab(\xb)),h_c(\betab(\xb))$, described in \Cref{sec:objective} subject to the constraints described in \Cref{sec:study}. The three constrained optimization problems are reformulated with an unconstrained quadratic penalty method and solved by using the manifold sampling trust-region method~\cite{Khan2017}. \Cref{fig:three_beta_plots} shows the resulting beta profiles for the optima under each objective. 
    
\begin{figure}[tb]
\centering
\includegraphics[width=0.8\linewidth]{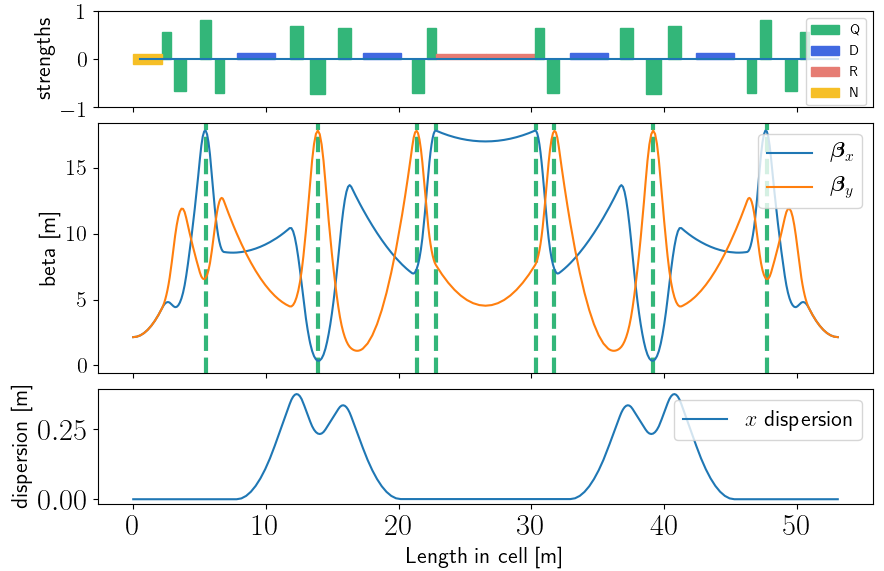}
\caption{(top) Visualization of the placement and strengths of the elements of the design $\xb_b$. Also shown are the corresponding beta profile (middle) and dispersion profile (bottom).
\label{fig:preflattice}}
\end{figure}

We find that optimization under the three objectives given the same initial point yields solutions that have notably similar beta profiles and correspondingly similar element positions, lengths, and strengths. 
% SW: Above seems kind of inverted from the usual way. I might have chopped sentence after "profiles"
While the accelerator design experts did not prefer one set of solutions over the other for their physical properties, they did prefer the solutions from $h_b$ because of the physical interpretability of the $\max_j\{\beta_j(\xb)\}$ objective. For the remainder of the exposition, we focus on the objective $h_b$ and the corresponding minimizer $\xb_b$ found, which is illustrated in \Cref{fig:preflattice}.

% \begin{table*}[!h]
% \small
% \begin{tabular}{rrrrrrrrrrr}
% \multicolumn{1}{l}{}          & \multicolumn{2}{c}{$h_a$}   & \multicolumn{2}{c}{$h_b$}     & \multicolumn{2}{c}{$h_c$}          
% \\ \hline
% \multicolumn{1}{l}{\textbf{}} & \multicolumn{1}{c}{\begin{tabular}[c]{@{}c@{}}Value \\ \end{tabular}} & \multicolumn{1}{c}{Activity} & \multicolumn{1}{c}{\begin{tabular}[c]{@{}c@{}}Value \\ \end{tabular}} & \multicolumn{1}{c}{Activity} & \multicolumn{1}{c}{\begin{tabular}[c]{@{}c@{}}Value\\ \end{tabular}} & \multicolumn{1}{c}{Activity} 
% \\ \hline
% Initial  &57.19  &(130, 112) & 29.91  &310 & 36.54 &33 \\
% Optima a  &38.21  &(47, 110) &19.55  &308 &27.27 &177 \\
% Optima b  &\textbf{34.61}  &(179, 111) & \textbf{17.31}  &309 &\textbf{25.32} &177 \\
% Optima c  &44.74  &(47, 111) &25.15  &309 &25.43 &160 \\
% \end{tabular}
% \caption{Optima for the objectives $h_a,h_b,h_c$. Optima a, indicates the optima for $h_a$. The activity is the index of the $\beta$-functions that is defines the objective.}
% \label{table:optima_table}
% \end{table*}

\subsection{A preferred lattice}
\label{sec:results-single_config}

\begin{figure}[tb]
\centering
\includegraphics[width=0.8\linewidth]{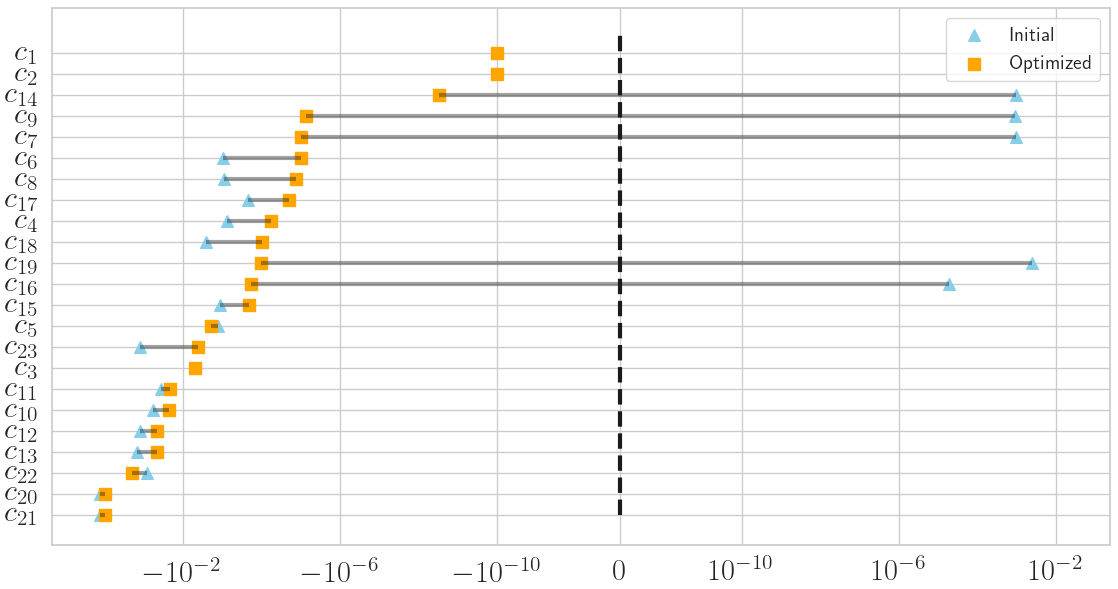}
\caption{Value of $c_i(\xb)$ for the initial lattice configuration $\xb=\eb_0$ (blue) and the optima $\xb=\xb_b$ (orange). Points left of the dotted black line indicate constraint satisfaction.}
\label{fig:lollipop_initial_optima}
\end{figure}

The beta profile for $\xb_b$ significantly improves on the beta profile of the initial lattice, reducing the max $\betab$ value from 29.91 to 17.80, a $40\%$ reduction. Furthermore, as shown in \Cref{fig:lollipop_initial_optima}, the $\xb_b$ lattice satisfies the constraints, whereas the initial lattice does not. Although several constraints are nearly active, this is in part due to the tight tolerances on the two-sided constraints, which imply that a small change in $\xb$ can easily change the constraint activities.
%Is this borne out in the later subsection?

\begin{table}[tbh]
    \centering
    \begin{tabular}{|l|r|r|}
        \hline
        \textbf{Parameter} & \textbf{Value} & \textbf{Units} \\
        \hline
        Particles & 1,000,000 & \\
         Normalized emittance & 24.0 & mm-mr\\
         Beta functions $x$ and $y$ & 2.15 & m \\
         Transverse RMS beam size & 2.34 & mm \\
         Width $\Delta p/p$ & 2.5\% & uniform\\
         \hline
    \end{tabular}
    \caption{Beam parameters for optimized lattice ($\xb_b$) simulation.}
    \label{tab:beam-param}
\end{table}

While the low beta functions and feasibility of our solution indicate that we have satisfied the design criteria, the true value of the candidate design is still unclear since the simulation used in this study  captures only linear optics. A full validation of our model and candidate solution leading up to a full project plan would require further investigation, particularly including a collective particle simulation.

\begin{figure}
    \centering
   \includegraphics[width=0.8\linewidth]{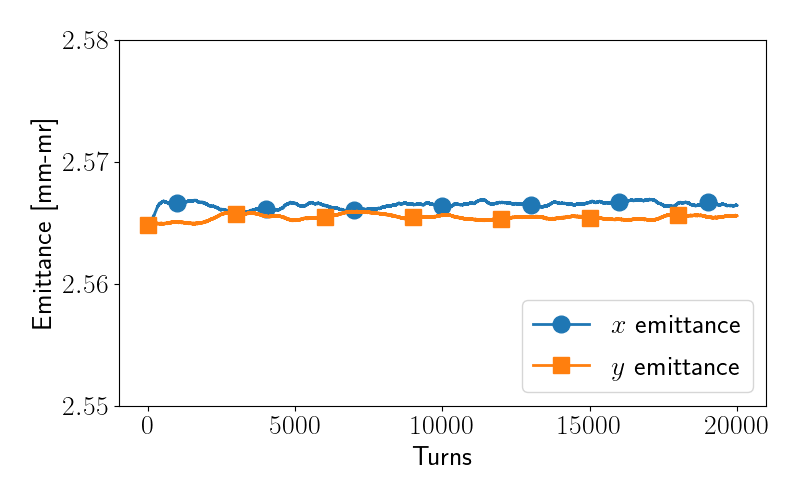}
    \caption{The $x$ and $y$ emittances of a representative beam distribution simulated using $\xb_b$. The beamline is very stable: even after 20,000 turns, the emittances have a relative change less than  0.1\%.}
    \label{fig:emittances}
\end{figure}

To satisfy ourselves that the optimized lattice would be able to propagate particles properly, we performed particle-tracking simulations using the $\xb_b$ lattice.
The beam parameters are shown in \Cref{tab:beam-param}.
Over the course of 20,000 simulated turns, there was no particle loss.
As seen in~\Cref{fig:emittances}, the emittance was stable within 0.1\% with no oscillations or long-term growth.

In applications, the constraints on the lattice optics should also be robust to small perturbations in the quadrupole strength. Perturbing the quadrupole strength simulates uncontrolled variation in magnet construction and power supplies. Any variation in dipole strength can be corrected with greater precision and has only indirect effect on the accelerator optics.

\begin{figure}[tb]
\centering
\includegraphics[width=0.8\linewidth]{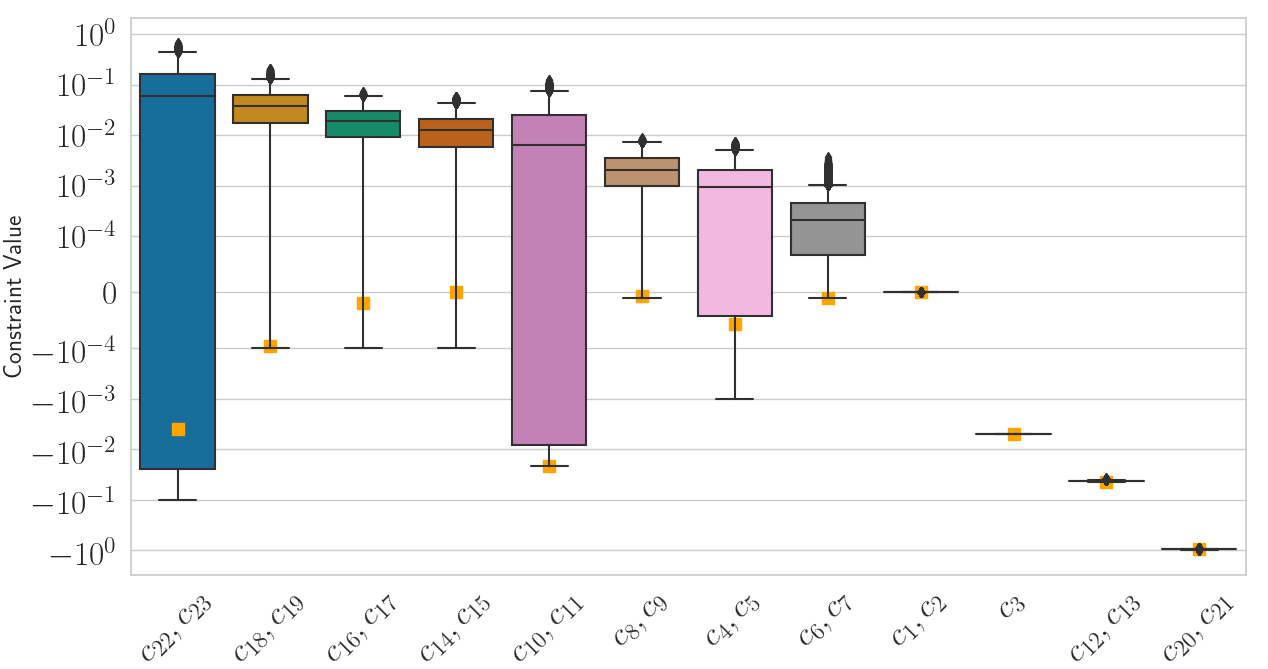}
\caption{Empirical distribution of constraint values $c_i(\xb_b + \mathbf{I}_q \bm{\xi})$ around $\xb_b$ under $\bm{\xi}\sim U([-10^{-3},10^{-3}]^{\nq})$ perturbations to the $\nq=9$ quadrupole strengths. When applicable, the maximum is taken over constraint pairs (upper and lower bounds) to show the worse of the two constraint violations. Orange squares indicate the constraint value (or maximum value for constraint pairs) at $\xb_b$. The distribution of constraint values under the random perturbations is robust since they are within a tolerable range to be corrected in tuning.}
\label{fig:perturbation_stability}
\end{figure}

\Cref{fig:perturbation_stability} shows that the $\xb_b$ lattice is indeed robust to uniformly distributed perturbations $U([-10^{-3},10^{-3}]^{\nq})$ to the $\nq$ quadrupole strengths. These perturbations reflect nontrivial changes of roughly 0.08\% relative to the original quadrupole strength bounds. The chromaticity matching criteria ($c_{22},c_{23}$) and the ratio of the beta functions in the center of the NL insert ($c_{10},c_{11}$) show the most significant variation in the constraint values and are easily restored under tuning. Those constraints, along with the dispersion-matching constraints ($c_{4}, \ldots, c_{9}$), exhibit perturbative behavior that naturally spans the constrained values. The betatron tune-matching ($c_{18},c_{19}$) and T-insert phase advances ($c_{14}, \ldots, c_{17}$) are Danilov--Nagaitsev integrability criteria that appeared to be violated under nearly any small random perturbation. Fortunately, those parameters can be measured with great precision in application and can be restored in dedicated tuning. The beam stability criterion ($c_{8},c_{9}$), the momentum-compaction factor ($c_{3}$), the beta matching at the edge of the NL insert ($c_{12},c_{13}$), and the phase advance over the NL insert ($c_{20},c_{21}$) are all insensitive to the small perturbations. 

The robustness of the lattice design $\xb_b$ may in part be due to the stringent simulation-based constraints required at the solution: we observe that points within the feasible region tend to be robust to large perturbations, while points violating  simulation-based constraints tend to demonstrate an instability in their constraint values when perturbed. 
%cite IOTA for tune-matching here?

\subsection{Sensitivity and tradeoff analysis}
\label{sec:results-tradeoffs}
The role of constraints in this problem cannot be overstated: 
all of the optima found are constrained optima in that the objective value would continue to decrease if the active constraints were removed. This begs the question of how constraint parameters, such as the right-hand side $\epsb$, play a role in selection of solutions. We explore this question through the KKT conditions at our constrained optima. The KKT conditions give rise to Lagrange multipliers, which quantify the improvement in objective value when active constraints are relaxed. 
We numerically solve \cref{eq:solve_for_lagrange} to compute the Lagrange multipliers. 

In order to be physically meaningful, a design must only be resolved so that the positions, lengths, and strengths are prescribed to a precision of roughly $10^{-3}$ or $10^{-4}$.
%units of the local optima. 
Element attributes cannot practically be tuned any tighter than this in a particle accelerator, so optimizing further is not meaningful from this perspective. However, in order to validate the KKT conditions and compute the Lagrange multipliers, solutions must be resolved to greater precision.

\begin{figure}[tb]
\centering
\includegraphics[width=0.8\linewidth]{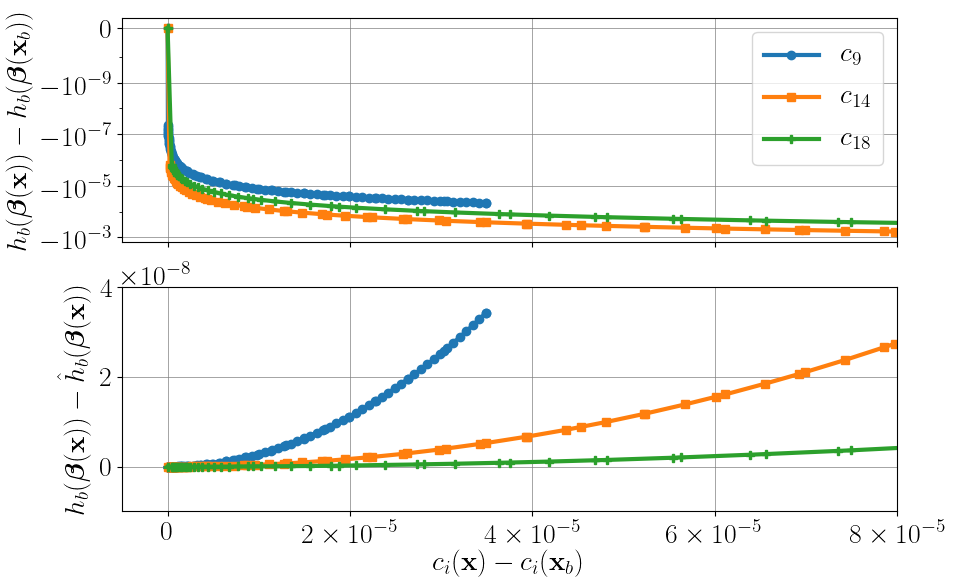}
\caption{(Top) Decrease of the optimal objective value as active constraints at $\xb_b$ are relaxed. (Bottom) The minuscule difference between $h_b$ and its predicted value using the linearization \cref{eq:predicted_decrease_linearization} while moving in the direction $\nabla c_i(\xb)$ for the active constraints. The high order of accuracy implies that the Lagrange multipliers capture the objective-constraint tradeoff tightly.}
\label{fig:pareto}
\end{figure}

The Lagrange multipliers for active constraints at  constrained optima give a prediction of how relaxation of these active constraints (i.e., increasing $\eps_i$ in \cref{eq:epsfeasible}) improves the objective value. For instance, at the optima $\xb_b$, where only simulation-based constraints indexed by $\mathcal{A}_{\rm sc}$ are active, the change in objective value along a direction $\db$ can be predicted locally with the linearization
\begin{equation}
    h_b(\betab(\xb)) \approx \hat{h}_b(\betab(\xb)) := h_b(\betab(\xb))  -\sum_{i\in\mathcal{A}_{\rm sc}(\xb_b)}\lambda_i^s\db^T\nabla c_i(\xb_b).
    \label{eq:predicted_decrease_linearization}
\end{equation}
While this shows that relaxing active constraints can improve the objective value, it also shows that relaxing inactive constraints will have no effect on the solution. In \Cref{fig:pareto} we see the actual and predicted decrease of the objective as we relax active constraints. 

\Cref{fig:pareto} shows that a small relaxation of constraint $c_{9}$ (the constraint on dispersion in the nonlinear section) yields a correspondingly small improvement in the objective. This suggests that a redesign of the achromatic bend may allow for smaller peak beta functions. Similarly, a small relaxation of constraint $c_{14}$ (phase advance over the T insert) or $c_{18}$ (matching the betatron tunes) also has a small impact on the objective. Because many of the Danilov--Nagaitsev integrability constraints are not active, the effect on the optimal solution of relaxing the complete set ($c_{6},\dots,c_{23}$) remains small. Furthermore, none of the bound or linear constraints are active at $\xb_b$, and so relaxing those will similarly not affect the solution.

In addition to interpreting the dependence of the optimal objective value on the active constraints, we can interpret the dependence of the optimal objective value on the active beta functions (i.e., the $\beta_j$ for which $h_b(\betab(\xb_b)) = \beta_j(\xb_b)$). At the optima $\xb_b$, the three beta functions $\beta_{x,48}$, $\beta_{x,180}$, $\beta_{y,112}$ are all essentially active, achieving values within $10^{-4}$ of the objective value $h_b(\betab(\xb_b))$. Mathematically, we can relax the dependence of the objective $h_b$ on these three beta function values by considering the alternative objective $h_b\left(\betab(\xb) - \sum_{j\in J^*(\xb_b)}\epsilon \eb_j\right)$,  where $\eb_j$ is the $j$th column of the identity matrix and $\epsilon>0$ is a small parameter to make the active beta function values inactive. The goal of relaxing active beta function values is to allow the objective to have an increased value at the beta function's associated element in order to have an improved performance along the remainder of the cell. 

\begin{figure}[tb]
\centering
\includegraphics[width=0.8\linewidth]{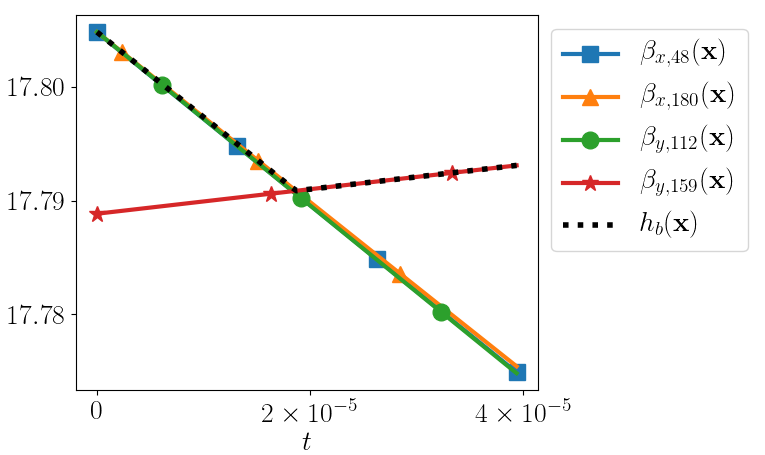}
\caption{Decrease in the objective $h_b$ and the three essentially active beta functions $\beta_{x,48},\beta_{x,180},\beta_{y,112}$ along the line segment $\xb_b + t\mathbf{d}^*$ until $\beta_{y,159}$ becomes active. At $\xb_b$, relaxing $\beta_{x,48},\beta_{x,180},\beta_{y,112}$ by roughly $\epsilon = 0.015$ would result in a change in activities, making $\beta_{y,159}$ the sole active index.}
\label{fig:objective_relaxation}
\end{figure}

\Cref{fig:objective_relaxation} shows the value of $h_b$, the three essentially active beta functions $\beta_{x,48}$, $\beta_{x,180}$, $\beta_{y,112}$, and $\beta_{y,159}$ along the line segment $\xb_b + t\mathbf{d}^*$, where $\mathbf{d}^*$ is the normalized solution to $\mbox{proj}\{ \zerob,\, \mbox{co}\{ \nabla \beta_{x,48}(\xb_b), \nabla \beta_{x,180}(\xb_b), \nabla \beta_{y,112}(\xb_b) \}\}$. At $\xb_b$, the three indices $\beta_{x,48}$, $\beta_{x,180}$, and $\beta_{y,112}$ would have to decrease by roughly $0.015$ in value before there would be a change in activities, making $\beta_{y,159}$ the sole active index. Furthermore under smaller relaxations the activities would not change, and so motion along the descent direction $\mathbf{d}^*$ would yield a brief decrease in $h_b$ before $\beta_{y,159}$ becomes active. This suggests a tradeoff between the peak value of the vertical beta function $\betab_{y}$ in the bending section ($\beta_{y,112}$) and the inner quadrupole doublet ($\beta_{y,159}$). 

The peak beta function can be allowed to be higher for a contiguous section of an accelerator lattice (such as the inner doublets and RF section) as a design choice in which the aperture is increased (which raises the cost and reduces the field strength of magnets in that section). As the peak vertical beta function in the bending section is reduced, the two peaks in the horizontal beta function $\betab_{x}$ within the bending section 
%(i.e., those associated with $\beta_{x,91}$ and $\beta_{x,130}$) 
may be increased without impacting $h_b$. Similarly, the peaks in the horizontal beta function within the bending section may be increased while the peak horizontal beta functions at the outer matching section ($\beta_{x,48}$) and inner quadrupole doublet ($\beta_{x,180}$) are reduced. 

%Rejected change to last section. Reducing the (maximum) field strength of the magnets is another adverse consequence of the increasing the aperture, not the purpose. The purpose is to allow the beta functions to be higher in that section (if the need to be). I.e. if the peak beta functions cannot be minimized uniformly across the accelerator, the next best thing would be to minimize them independently for different section.

\section{Conclusion}\label{sec:conc}
Through the development and solution of an optimization model we explored the design of an integrable rapid cycling synchrotron. The model leverages linear optics to rapidly simulate the properties of lattice designs. Challenges inherit to the model include its sizable dimension (32 decision variables, 55 linear constraints, and 23 simulation-based constraints); lack of derivative availability of the simulation-based quantities; and a nonsmooth, simulation-based objective function. By judicious navigation of the feasible region and exploitation of known compositions of the simulation-based quantities, application of a manifold sampling algorithm yielded solutions that not only verify (up to numerical tolerances) the viability of the integrable lattice design but also perform well. 
This study aims to serve as a foundation for further accelerator optimization studies by methodically formulating a model, taking advantage of problem structure, and studying key sensitivities. 

Additional improvements can be made toward a more comprehensive and general-purpose design of synchrotrons. 
For example, this optimization relied on specific optics-based lattice criteria that prior Synergia simulations of intense nonlinear beams have shown to be associated with achieving benchmarks in machine performance, such as the beam quality and the loss-limited beam intensity. An analysis complementary to our work could integrate nonlinear space-charge simulations into the optimization model directly to explore the relationship between machine parameters and those performance benchmarks. However, the outcome of such an optimization would necessarily be contingent on assumptions regarding the initial beam distribution, the allocation of machine errors, the availability of beam tuning, and the particle loss model.

The lattice optimization also takes place within the context of an overall lattice cell structure (the sequence of magnets as well as the number of periodic cells). In prior design work, this overall lattice cell structure was generated manually through a combination of domain knowledge and elementary operations (splitting, merging, transposing, and changing the number of periodic cells). If the lattice optimization presented here could be paired with an algorithm for selecting and rejecting cell structures, then a truly general accelerator optimization result could be generated. A key challenge in such a method is that it would require efficiently managing the vast majority of cell structures that are unstable, directly incompatible with constraints, or highly ineffective.

\section*{Acknowledgements}
We gratefully acknowledge the computing resources provided on Bebop, a high-performance computing cluster operated by the Laboratory Computing Resource Center at Argonne National Laboratory.


\bibliographystyle{spmpsci}      % mathematics and physical sciences
\bibliography{bibs/smw-bigrefs.bib,bibs/allbibs.bib,bibs/egs_bibs.bib}   

\begin{thebibliography}{10}
\providecommand{\url}[1]{{#1}}
\providecommand{\urlprefix}{URL }
\expandafter\ifx\csname urlstyle\endcsname\relax
  \providecommand{\doi}[1]{DOI~\discretionary{}{}{}#1}\else
  \providecommand{\doi}{DOI~\discretionary{}{}{}\begingroup
  \urlstyle{rm}\Url}\fi

\bibitem{DUNE}
Abi, R., et~al.: {Deep Underground Neutrino Experiment (DUNE), Far Detector
  Technical Design Report, Volume II: DUNE Physics}.
\newblock Fermilab, Batavia, FERMILAB-PUB-20-025-ND  (2020).
\newblock \urlprefix\url{https://arxiv.org/abs/2002.03005}

\bibitem{EldredSyphers}
Ainsworth, R., Dey, J., Eldred, J., Harnik, R., Jarvis, J., Johnson, D.E.,
  Kourbanis, I., Neuffer, D., Pozdeyev, E., Syphers, M.J., Valishev, A.,
  Yakovlev, V.P., Zwaska, R.: An upgrade path for the {Fermilab} accelerator
  complex.
\newblock Fermilab, Batavia, FERMILAB-TM-2754-AD-APC-PIP2-TD  (2021).
\newblock \urlprefix\url{https://arxiv.org/abs/2106.02133}

\bibitem{Synergia}
Amundson, J., Goldhaber, S., Lebrun, P., Lu, Q., Macridin, A., Michelotti, L.,
  Park, C.S., Spentzouris, P., Stern, E.: Synergia.
\newblock \urlprefix\url{https://synergia.fnal.gov/}

\bibitem{Amundson:2004qd}
Amundson, J.F., Spentzouris, P., Qiang, J., Ryne, R.: Synergia: an accelerator
  modeling tool with 3-{D} space charge.
\newblock Journal of Computational Physics \textbf{211}, 229--248 (2006).
\newblock \doi{10.1016/j.jcp.2005.05.024}

\bibitem{Antipov2017}
Antipov, S., Broemmelsiek, D., Bruhwiler, D., Edstrom, D., Harms, E., Lebedev,
  V., Leibfritz, J., Nagaitsev, S., Park, C., Piekarz, H., Piot, P., Prebys,
  E., Romanov, A., Ruan, J., Sen, T., Stancari, G., Thangaraj, C.,
  Thurman-Keup, R., Valishev, A., Shiltsev, V.: {IOTA} {(Integrable Optics Test
  Accelerator)}: Facility and experimental beam physics program.
\newblock Journal of Instrumentation \textbf{12}(03), T03002--T03002 (2017).
\newblock \doi{10.1088/1748-0221/12/03/t03002}

\bibitem{AudetHare2017}
Audet, C., Hare, W.L.: Derivative-Free and Blackbox Optimization.
\newblock Springer (2017).
\newblock \doi{10.1007/978-3-319-68913-5}

\bibitem{bagirov2014intrononsmooth}
Bagirov, A., Karmitsa, N., M{\"a}kel{\"a}, M.M.: Introduction to Nonsmooth
  Optimization: Theory, Practice and Software.
\newblock Springer (2014).
\newblock \doi{10.1007/978-3-319-08114-4}

\bibitem{bueler2000vinci}
B{\"u}eler, B., Enge, A.: Vinci (2000).
\newblock \urlprefix\url{https://www.math.u-bordeaux.fr/~aenge}

\bibitem{chen2018polytope}
Chen, Y., Dwivedi, R., Wainwright, M.J., Yu, B.: Fast {MCMC} sampling
  algorithms on polytopes.
\newblock The Journal of Machine Learning Research \textbf{19}(1), 2146--2231
  (2018).
\newblock \urlprefix\url{https://jmlr.org/papers/v19/18-158.html}

\bibitem{Danilov}
Danilov, V., Nagaitsev, S.: Nonlinear accelerator lattices with one and two
  analytic invariants.
\newblock Physical Review Accelerators and Beams \textbf{13}(8), 084002 (2010).
\newblock \doi{10.1103/PhysRevSTAB.13.084002}

\bibitem{Duris2020}
Duris, J., Kennedy, D., Hanuka, A., Shtalenkova, J., Edelen, A., Baxevanis, P.,
  Egger, A., Cope, T., McIntire, M., Ermon, S., Ratner, D.: Bayesian
  optimization of a free-electron laser.
\newblock Physical Review Letters \textbf{124}(12) (2020).
\newblock \doi{10.1103/physrevlett.124.124801}

\bibitem{Edelen2020}
Edelen, A., Neveu, N., Frey, M., Huber, Y., Mayes, C., Adelmann, A.: Machine
  learning for orders of magnitude speedup in multiobjective optimization of
  particle accelerator systems.
\newblock Physical Review Accelerators and Beams \textbf{23}(4) (2020).
\newblock \doi{10.1103/physrevaccelbeams.23.044601}

\bibitem{SyphersBook}
Edwards, D.A., Syphers, M.J.: An Introduction to the Physics of High Energy
  Accelerators.
\newblock Wiley-VCH (1993).
\newblock \doi{10.1002/9783527617272}

\bibitem{Eldred2020}
Eldred, J.: Novel approaches to high-power proton beams.
\newblock {PoS} \textbf{NuFact2019}, 055 (2020).
\newblock \doi{10.22323/1.369.0055}

\bibitem{EldredJINST}
Eldred, J., Lebedev, V., Valishev, A.: Rapid-cycling synchrotron for
  multi-megawatt proton facility at {Fermilab}.
\newblock {Journal} of Instrumentation \textbf{14}(07), P07021 (2019).
\newblock \doi{10.1088/1748-0221/14/07/P07021}

\bibitem{EldredIPAC17}
Eldred, J., Valishev, A.: Space-charge simulation of integrable rapid cycling
  synchrotron.
\newblock Proceedings of the 8th Int. Particle Accelerator Conf.  (2017).
\newblock \doi{10.18429/JACOW-IPAC2017-THPVA032}

\bibitem{EldredIPAC18}
Eldred, J., Valishev, A.: Simulation of integrable synchrotron with
  space-charge and chromatic tune-shifts.
\newblock Proceedings of the 9th Int. Particle Accelerator Conf.  (2018).
\newblock \doi{10.18429/JACOW-IPAC2018-TUPAF073}

\bibitem{Fletcher}
Fletcher, R.: Practical Methods of Optimization, second edn.
\newblock John Wiley \& Sons (1987).
\newblock \doi{10.1002/9781118723203}

\bibitem{Huang2018}
Huang, X.: Robust simplex algorithm for online optimization.
\newblock Physical Review Accelerators and Beams \textbf{21}(10), 104601
  (2018).
\newblock \doi{10.1103/PhysRevAccelBeams.21.104601}

\bibitem{HuangSafranek2014}
Huang, X., Safranek, J.: Nonlinear dynamics optimization with particle swarm
  and genetic algorithms for {SPEAR}3 emittance upgrade.
\newblock Nuclear Instruments and Methods in Physics Research Section A:
  Accelerators, Spectrometers, Detectors and Associated Equipment \textbf{757},
  48--53 (2014).
\newblock \doi{10.1016/j.nima.2014.04.078}

\bibitem{Khachiyan2008}
Khachiyan, L., Boros, E., Borys, K., Elbassioni, K., Gurvich, V.: Generating
  all vertices of a polyhedron is hard.
\newblock Discrete {\&} Computational Geometry \textbf{39}(1-3), 174--190
  (2008).
\newblock \doi{10.1007/s00454-008-9050-5}

\bibitem{Khan2017}
Khan, K.A.: Branch-locking {AD} techniques for nonsmooth composite functions
  and nonsmooth implicit functions.
\newblock Optimization Methods and Software \textbf{33}(4-6), 1127--1155
  (2017).
\newblock \doi{10.1080/10556788.2017.1341506}

\bibitem{KLW18}
Khan, K.A., Larson, J., Wild, S.M.: Manifold sampling for optimization of
  nonconvex functions that are piecewise linear compositions of smooth
  components.
\newblock {SIAM} Journal on Optimization \textbf{28}(4), 3001--3024 (2018).
\newblock \doi{10.1137/17m114741x}

\bibitem{LMW16}
Larson, J., Menickelly, M., Wild, S.M.: Manifold sampling for $\ell_1$
  nonconvex optimization.
\newblock SIAM Journal on Optimization \textbf{26}(4), 2540--2563 (2016).
\newblock \doi{10.1137/15M1042097}

\bibitem{LMW2019AN}
Larson, J., Menickelly, M., Wild, S.M.: Derivative-free optimization methods.
\newblock Acta Numerica \textbf{28}, 287--404 (2019).
\newblock \doi{10.1017/s0962492919000060}

\bibitem{Larson2020}
Larson, J., Menickelly, M., Zhou, B.: Manifold sampling for optimizing
  nonsmooth nonconvex compositions.
\newblock SIAM Journal on Optimization  (2021).
\newblock \urlprefix\url{https://arxiv.org/abs/2011.01283}.
\newblock To appear

\bibitem{taxonomy15}
{Le Digabel}, S., Wild, S.M.: A taxonomy of constraints in black-box
  simulation-based optimization.
\newblock Preprint ANL/MCS-P5350-0515, Argonne National Laboratory, Mathematics
  and Computer Science Division (2015-01).
\newblock \urlprefix\url{http://www.mcs.anl.gov/papers/P5350-0515.pdf}

\bibitem{PIP2}
Lebedev, V., et~al.: The {PIP-II} conceptual design report.
\newblock Fermilab, Batavia, FERMILAB-TM-2649-AD-APC  (2017).
\newblock \urlprefix\url{https://pxie.fnal.gov/PIP-II_CDR/PIP-II_CDR_v.0.3.pdf}

\bibitem{LeeBook}
Lee, S.Y.: Accelerator Physics, third edn.
\newblock World Scientific (2011).
\newblock \doi{10.1142/8335}

\bibitem{Li2018}
Li, Y., Cheng, W., Yu, L.H., Rainer, R.: Genetic algorithm enhanced by machine
  learning in dynamic aperture optimization.
\newblock Physical Review Accelerators and Beams \textbf{21}(5), 054601 (2018).
\newblock \doi{10.1103/PhysRevAccelBeams.21.054601}

\bibitem{Macridin:2015vua}
Macridin, A., Burov, A., Stern, E., Amundson, J., Spentzouris, P.: Simulation
  of transverse modes with their intrinsic landau damping for bunched beams in
  the presence of space charge.
\newblock Physical Review Accelerators and Beams \textbf{18}(7), 074401 (2015).
\newblock \doi{10.1103/PhysRevSTAB.18.074401}

\bibitem{Martin2019}
Martin, L.K., Kelliher, D.J., Sheehy, S.L.: Can a {Paul} ion trap be used to
  investigate nonlinear quasi-integrable optics?
\newblock Journal of Physics: Conference Series \textbf{1350}, 012132 (2019).
\newblock \doi{10.1088/1742-6596/1350/1/012132}

\bibitem{mete2012pattern}
Mete, H.O., Zabinsky, Z.B.: Pattern hit-and-run for sampling efficiently on
  polytopes.
\newblock Operations Research Letters \textbf{40}(1), 6--11 (2012).
\newblock \doi{10.1016/j.orl.2011.11.002}

\bibitem{Michelotti2006}
Michelotti, L., Ostiguy, J.F.: {CHEF: A} framework for accelerator optics and
  simulation.
\newblock In: Proceedings of 9th International Computational Accelerator
  Physics Conference, pp. 2--6 (2006)

\bibitem{MintyBook}
Minty, M.G., Zimmermann, F.: Measurement and Control of Charged Particle Beams.
\newblock Springer (2003).
\newblock \doi{10.1007/978-3-662-08581-3}

\bibitem{Nagaitsev}
Nagaitsev, S., Lebedev, V.: A cost-effective rapid-cycling synchrotron.
\newblock Reviews of Accelerator Science and Technology \textbf{10}(01),
  245--266 (2019).
\newblock \doi{10.1142/s1793626819300135}

\bibitem{NelderMead}
Nelder, J.A., Mead, R.: A simplex method for function minimization.
\newblock The Computer Journal \textbf{7}(4), 308--313 (1965).
\newblock \doi{10.1093/comjnl/7.4.308}

\bibitem{Neveu2019}
Neveu, N., Hudson, S., Larson, J., Spentzouris, L.: Comparison of model-based
  and heuristic optimization algorithms applied to photoinjectors using
  {libEnsemble}.
\newblock In: Proceedings of the 13th International Computational Accelerator
  Physics Conference, pp. 22--24 (2019).
\newblock \doi{10.18429/JACoW-ICAP2018-SAPAF03}

\bibitem{Neveu2017}
Neveu, N., Larson, J., Power, J.G., Spentzouris, L.: Photoinjector optimization
  using a derivative-free, model-based trust-region algorithm for the {Argonne}
  {Wakefield} {Accelerator}.
\newblock Journal of Physics: Conference Series \textbf{874}, 012062 (2017).
\newblock \doi{10.1088/1742-6596/874/1/012062}

\bibitem{nocedal2006no}
Nocedal, J., Wright, S.J.: Numerical Optimization, second edn.
\newblock Springer (2006).
\newblock \doi{10.1007/978-0-387-40065-5}

\bibitem{Pang2014}
Pang, X., Rybarcyk, L.: Multi-objective particle swarm and genetic algorithm
  for the optimization of the {LANSCE} linac operation.
\newblock Nuclear Instruments and Methods in Physics Research Section A:
  Accelerators, Spectrometers, Detectors and Associated Equipment \textbf{741},
  124--129 (2014).
\newblock \doi{10.1016/j.nima.2013.12.042}

\bibitem{Powell1994}
Powell, M.J.D.: A direct search optimization method that models the objective
  and constraint functions by linear interpolation.
\newblock In: S.~Gomez, J.P. Hennart (eds.) Advances in Optimization and
  Numerical Analysis, \emph{Mathematics and its Applications}, vol. 275, pp.
  51--67. Springer (1994).
\newblock \doi{10.1007/978-94-015-8330-5_4}

\bibitem{pdfo}
Ragonneau, T.M., Zhang, Z.: {PDFO}: Cross-platform interfaces for {Powell's}
  derivative-free optimization solvers (version 1.0).
\newblock \doi{10.5281/zenodo.3887569}

\bibitem{Roussel2021}
Roussel, R., Hanuka, A., Edelen, A.: Multiobjective {Bayesian} optimization for
  online accelerator tuning.
\newblock Phys. Rev. Accel. Beams \textbf{24}, 062801 (2021).
\newblock \doi{10.1103/PhysRevAccelBeams.24.062801}

\bibitem{Ruisard2019}
Ruisard, K., Komkov, H.B., Beaudoin, B., Haber, I., Matthew, D., Koeth, T.:
  Single-invariant nonlinear optics for a small electron recirculator.
\newblock Physical Review Accelerators and Beams \textbf{22}(4), 041601 (2019).
\newblock \doi{10.1103/PhysRevAccelBeams.22.041601}

\bibitem{Scheinker2014}
Scheinker, A., Pang, X., Rybarcyk, L.: Model-independent particle accelerator
  tuning.
\newblock Physical Review Accelerators and Beams \textbf{16}(10), 102803
  (2014).
\newblock \doi{10.1103/PhysRevSTAB.16.102803}

\bibitem{Borland2005}
Shang, H., Borland, M.: A parallel simplex optimizer and its application to
  high-brightness storage ring design.
\newblock In: Proceedings of the 2005 Particle Accelerator Conference, pp.
  4230--4232. IEEE (2005).
\newblock \doi{10.1109/pac.2005.1591774}

\bibitem{Shiltsev2017}
Shiltsev, V.: Fermilab proton accelerator complex status and improvement plans.
\newblock Modern Physics Letters A \textbf{32}(16) (2017).
\newblock \doi{10.1142/S0217732317300129}

\bibitem{Sun2017}
Sun, Y.: Multi-objective online optimization of beam lifetime at {APS}.
\newblock In: Proceedings of North American Particle Accelerator Conference, 3,
  pp. 913--915 (2017).
\newblock \doi{https://doi.org/10.18429/JACoW-NAPAC2016-WEPOB12}

\bibitem{Valishev2020}
Valishev, A.: {Research at FAST/IOTA: Strategy and priorities}.
\newblock IOTA Collaboration Meeting (2020).
\newblock
  \urlprefix\url{https://indico.fnal.gov/event/43231/contributions/187342/attachments/129553/157411/2020-06-15_Strategy_CollaborationMeeting.pdf}

\bibitem{Webb2020}
Webb, S., Cook, N., Eldred, J.: Averaged invariants in storage rings with
  synchrotron motion.
\newblock Journal of Instrumentation \textbf{15}(12), 12032 (2020).
\newblock \doi{10.1088/1748-0221/15/12/p12032}

\bibitem{Webb2012}
Webb, S.D., Bruhwiler, D.L., Abell, D.T., Sishlo, A., Danilov, V., Nagaitsev,
  S., Valishev, A., Danilov, K., Cary, J.R.: Effects of nonlinear decoherence
  on halo formation.
\newblock Tech. Rep. 1205.7083, ArXiv (2012).
\newblock \urlprefix\url{https://arxiv.org/abs/1205.7083}

\bibitem{Webb2015}
Webb, S.D., Bruhwiler, D.L., Valishev, A., Nagaitsev, S.N., Danilov, V.V.:
  Chromatic and dispersive effects in nonlinear integrable optics.
\newblock Tech. Rep. 1504.05981, ArXiv (2015).
\newblock \urlprefix\url{https://arxiv.org/abs/1504.05981}

\bibitem{Wei2003}
Wei, J.: Synchrotrons and accumulators for high-intensity proton beams.
\newblock Reviews of Modern Physics \textbf{75}(1383), 1383--1432 (2003).
\newblock \doi{10.1103/RevModPhys.75.1383}

\bibitem{Yang2011}
Yang, L., Li, Y., Guo, W., Krinsky, S.: Multiobjective optimization of dynamic
  aperture.
\newblock Physical Review Accelerators and Beams \textbf{14}(5), 054001 (2011).
\newblock \doi{10.1103/PhysRevSTAB.14.054001}

\end{thebibliography}

% SW: This should be deleted in the version that appears, but we will need to show it for purposes of clearing the manuscript at Argonne and it's usually not a bad thing to show in preprint versions.
%\clearpage
\vspace{3em}

\small
\framebox{\parbox{\linewidth}{
%The submitted manuscript has been created by UChicago Argonne, LLC, Operator of Argonne National Laboratory (``Argonne''). Argonne, a U.S.\ Department of Energy Office of Science laboratory, is operated under Contract No.\ DE-AC02-06CH11357. 
The U.S.\ Government retains for itself, and others acting on its behalf, a paid-up nonexclusive, irrevocable worldwide license in said article to reproduce, prepare derivative works, distribute copies to the public, and perform publicly and display publicly, by or on behalf of the Government.  The Department of Energy will provide public access to these results of federally sponsored research in accordance with the DOE Public Access Plan. http://energy.gov/downloads/doe-public-access-plan.}}

\end{document}